\documentclass[aos,MSNbibl,nameyear,dvips]{arximspdf}
\usepackage{graphicx}
\doi{10.1214/13-AOS1182} 
\volume{41}
\issue{6}
\pubyear{2013}
\firstpage{2994}
\lastpage{3021}

\makeatletter

\def\cal{\mathcal}
\newcommand{\ve}[0]{\mathbf{e}}
\newcommand{\vg}[0]{\mathbf{g}}
\newcommand{\vq}[0]{\mathbf{q}}
\newcommand{\vx}[0]{\mathbf{x}}
\newcommand{\vz}[0]{\mathbf{z}}
\newcommand{\diag}[0]{\operatorname{diag}}
\newcommand{\tr}[0]{\operatorname{tr}}
\newcommand{\calF}[0]{\mathcal{F}}
\newcommand{\calH}[0]{\mathcal{H}}
\newcommand{\E}[0]{\mathbb{E}}
\newcommand{\BL}[0]{\mathrm{BL}}

\newcommand{\Prob}[0]{\mathbb{P}}
\newcommand{\I}[0]{\mathbb{I}}
\renewcommand{\vec}[0]{\operatorname{vec}}

\newtheorem{thmm}{Theorem}[section]

\newtheorem{cor}[thmm]{Corollary}

\newproclaim{defn}{Definition}[section]
\newproclaim{conj}{Conjecture}[section]
\newproclaim{exmp}{Example}[section]

\newproclaim{rmk}{Remark}
\newproclaim{note}{Note}
\newproclaim{case}{Case}
\makeatother

\begin{document}
\begin{frontmatter}

\title{Covariance and precision matrix estimation for high-dimensional time series}
\runtitle{High-dimensional covariance estimation for time series}

\begin{aug}
\author[a]{\fnms{Xiaohui} \snm{Chen}\corref{}\ead[label=e1]{xhchen@illinois.edu}},
\author[b]{\fnms{Mengyu} \snm{Xu}\ead[label=e2]{mengyu@galton.uchicago.edu}}
\and
\author[b]{\fnms{Wei Biao} \snm{Wu}\ead[label=e3]{wbwu@galton.uchicago.edu}}
\runauthor{X. Chen, M. Xu and W.~B. Wu}
\affiliation{University of Illinois at Urbana-Champaign, University of
Chicago and
University of Chicago}
\address[a]{X. Chen\\
Department of Statistics\\
University of Illinois at Urbana-Champaign\\
725 S. Wright Street\\
Champaign, Illinois 61820\\
USA\\
\printead{e1}}

\address[b]{M. Xu\\
W.~B. Wu\\
Department of Statistics\\
University of Chicago\\
5734 S. University Avenue\\
Chicago, Illinois 60637\\
USA\\
\printead{e2} \\
\phantom{E-mail:\ }\printead*{e3}}
\end{aug}

\received{\smonth{4} \syear{2013}}
\revised{\smonth{10} \syear{2013}}

%
\begin{abstract}
We consider estimation of covariance matrices and their inverses
(a.k.a. precision matrices) for high-dimensional stationary and
locally stationary time series. In the latter case the covariance
matrices evolve smoothly in time, thus forming a covariance matrix
function. Using the functional dependence measure of Wu
[\textit{Proc. Natl. Acad. Sci. USA}
\textbf{102} (2005) 14150--14154 (electronic)],
we obtain the rate of convergence for the
thresholded estimate and illustrate how the dependence affects the
rate of convergence. Asymptotic properties are also obtained for
the precision matrix estimate which is based on the graphical
Lasso principle. Our theory substantially generalizes earlier ones
by allowing dependence, by allowing nonstationarity and by
relaxing the associated moment conditions.
\end{abstract}

%

\end{frontmatter}

\section{Introduction}
\label{sec:introduction}

Estimation of covariance matrices and their inverses (a.k.a.
precision matrices) is of fundamental importance in almost every
aspect of statistics, ranging from the principal component
analysis [\citet{johnstonelu2009}], graphical modeling
[\citet{meinshausenbuhlmann2006,ravikumarwainwrightraskuttiyu2008a,yuan2010a}],
classification
based on the linear or quadratic discriminant analysis
[\citet{bickellevina2004}], and real-world applications such as
portfolio selection [\citet{ledoitwolf2003a,talih2003}] and
wireless communication [\citet{guerci1999a,ward1994a,listociawang2003a,abrahamssonselenstoica2007a}]. Suppose we have
$n$ temporally observed $p$-dimensional vectors $(\vz_i)_{i=1}^n$,
with $\vz_i$ having mean zero and covariance matrix $\Sigma_i =
\E(\vz_i \vz_i^\top)$ whose dimension is $p \times p$. Our goal is
to estimate the covariance matrices $\Sigma_i$ and their inverses
$\Omega_i = \Sigma_i^{-1}$ based on the data matrix $Z_{p \times
n} = (\vz_1, \ldots, \vz_n)$. In the classical situation where $p$
is fixed, $n \to\infty$ and $\vz_i$ are mean zero independent and
identically distributed (i.i.d.) random vectors, it is well known
that the sample covariance matrix
%
\begin{equation}
\label{eq:samplecovariancematrix} \hat\Sigma_n = n^{-1} \sum
_{i=1}^n \vz_i\vz_i^\top
\end{equation}
is a consistent and well behaved estimator of $\Sigma$, and $\hat
\Omega_n = \hat\Sigma_n^{-1}$ is a natural and good estimator of
$\Omega$. See \citet{MR0091588} for a detailed account. However,
when the dimensionality $p$ grows with $n$, random matrix theory
asserts that $\hat\Sigma_n$ is no longer a consistent estimate of
$\Sigma$ in the sense that its eigenvalues do not converge to
those of $\Sigma$; see, for example, the Mar{\v{c}}enko--Pastur law
[\citet{marchenkopastur1967}] or the Tracy--Widom law
[\citet{MR1863961}]. Moreover, it is clear that $\hat\Omega_n$ is
not defined when $\hat\Sigma_n$ is not invertible in the
high-dimensional case with $p > n$.

During the last decade, various special cases of the above
covariance matrix estimation problem have been studied. In most of
the previous papers it is assumed that the vectors $\vz_1, \ldots,
\vz_n$ are i.i.d. and thus the covariance matrix $\Sigma_i \equiv
\Sigma$ is time-invariant. See, for example,
\citeauthor{bickellevina2008a} (\citeyear{bickellevina2008a,bickellevina2008b}), \citet{caizhangzhou2010a},
\citeauthor{caizhou2011a} (\citeyear{caizhou2011a,caizhou2011b}), where consistency and rates of
convergence are established for various regularized (banded,
tapered or thresholded) estimates of covariance matrices and their
inverses. As an alternative regularized estimate for sparse
precision matrix, one can adopt the Lasso-type entry-wise 1-norm
penalized likelihood approach; see
\citet{rothmanbickellevinazhu2008a,friedmanhastietibshirani2008a,banerjeeelghaouidaspremont2008a,ravikumarwainwrightraskuttiyu2008a,fanfengwu2009a}. Other
estimates include the Cholesky decomposition based method
[\citet{wupourahmadi2003a,huangliupourahmadiliu2006a}],
neighborhood selection for sparse graphical models
[\citet{liuluo2012,yuan2010a,meinshausenbuhlmann2006}], regularized
likelihood approach [\citet{lamfan2009a,fanfengwu2009a}] and the
sparse matrix transform [\citet{caobachegabouman2011}].
\citet{Xiao2012} considered covariance matrix estimation for
univariate stationary processes.

The assumption that $\vz_1, \ldots, \vz_n$ are i.i.d. is quite
restrictive for situations that involve temporally observed data.
In \citet{zhoulaffertywasserman2010a} and \citet{kolarxing2011a} the
authors considered time-varying Gaussian graphical models where
the sampling distribution can change smoothly over time. However,
they assume that the underlying random vectors are independent.
Using nonparametric smoothing techniques, they estimate the
time-vary covariance matrices in terms of covariance matrix
functions. Their asymptotic theory critically depends on the \textit
{independence} assumption.

The importance of estimating covariance matrices for dependent and
nonstationary processes has been increasingly seen across a wide
variety of research areas. In modeling spatial--temporal data,
\citet{wiklehooten2010a} proposed quadratic nonlinear dynamic
models to accommodate the interactions between the processes which
are useful for characterizing dynamic processes in geophysics
[\citet{Kondrashov2005a}]. \citet{zhengchenblasch2007} considered
non-Gaussian clutter and noise processes in space--time adaptive
processing, where the space--time covariance matrix is important
for detecting airborne moving targets in the nonstationary
clutter environment [\citet{ward1994a,guerci1999a}]. In finance,
\citet{jacquierpolsonrossi2004} considered multivariate stochastic
volatility models parametrized by time-varying covariance matrices
with heavy tails and correlated errors. \citet{talih2003}
investigated the Markowitz portfolio selection problem for optimal
returns of a large number of stocks with hidden and heterogeneous
Gaussian graphical model structures. In essence, those real-world
problems pose a number of challenges: (i) nonlinear dynamics of
data generating systems, (ii)~temporally dependent and
nonstationary observations, (iii) high-dimensionality of the
parameter space and (iv) non-Gaussian distributions. Therefore, the
combination of more flexible nonlinear and nonstationary
components in the models and regularized covariance matrix
estimation are essential to perform related statistical inference.

In contrast to the longstanding progresses and extensive research
that have been made in terms of heuristics and methodology,
theoretical work on estimation of covariance matrices based on
high-dimensional time series data is largely untouched. In this
paper we shall substantially relax the i.i.d. assumption by
establishing an asymptotic theory that can have a wide range of
applicability. We shall deal with the estimation of covariance and
precision matrices for high-dimensional stationary processes in
Sections~\ref{sec:stationary} and \ref{sec:precision_statproc},
respectively. Section~\ref{sec:stationary} provides a rate of
convergence for the thresholded estimator, and Section~\ref{sec:precision_statproc} concerns the graphical Lasso
estimator for precision matrices. For locally stationary
processes, an important class of nonstationary processes, we
shall study in Section~\ref{sec:covariancematrixestiamtion_nonstatproc} the estimation of
time-varying covariance and precision matrices. This
generalization allows us to consider time-varying covariance and
precision matrix estimation under temporal dependence; hence our
results significantly extend previous ones by
\citet{zhoulaffertywasserman2010a} and \citet{kolarxing2011a}.
Furthermore, by assuming a mild moment condition on the underlying
processes, we can relax the multivariate Gaussian assumption that
was imposed in \citet{zhoulaffertywasserman2010a}
and~\citet{kolarxing2011a} [and also by \citeauthor{bickellevina2008a}
(\citeyear{bickellevina2008a,bickellevina2008b})
in the i.i.d. setting]. Specifically, we shall
show that, thresholding on the kernel smoothed sample covariance
matrices, estimators based on the localized graphical Lasso
procedure are consistent estimators for time-varying covariance
and precision matrices.

To deal with temporal dependence, we shall use the functional
dependence measure of \citet{MR2172215}. With the latter, we are
able to obtain explicit rates of convergence for the thresholded
covariance matrix estimates and illustrate how the dependence
affects the rates. In particular, we show that, based on the
moment condition of the underlying process, there exists a
threshold value. If the dependence of the process does not exceed
that threshold, then the rates of convergence will be the same as
those obtained under independence. On the other hand, if the
dependence is stronger, then the rates of convergence will depend
on the dependence. This phase transition phenomenon is of
independent interest.

We now introduce some notation. We shall use $C, C_1, C_2, \ldots $
to denote positive constants whose values may differ from place to
place. Those constants are independent of the sample size $n$ and
the dimension $p$. For some quantities $a$ and $b$, which may
depend on $n$ and $p$, we write $a \lesssim b$ if $a \le C b$
holds for some constant $C$ that is independent of $n$ and $p$ and
$a \asymp b$ if there exists a constant $0 < C < \infty$ such that
$C \le\liminf b/a \le\limsup b/a \le C^{-1}$. We use $x \wedge y
= \min(x,y)$ and $x \vee y = \max(x,y)$. For a vector $\vx\in
\mathbb{R}^p$, we write $|\vx|= (\sum_{j=1}^p x_j^2)^{1/2}$ and
for a matrix\vspace*{1pt} $\Sigma$, $|\Sigma|_1 = \sum_{j,k}|\sigma_{j k}|$,
$|\Sigma|_\infty= \max_{j,k} |\sigma_{jk}|$, $|\Sigma|_F =
(\sum_{j, k}\sigma_{j k}^2)^{1/2}$ and $\rho(\Sigma) = \max
\{|\Sigma\vx| \dvtx |\vx|=1\}$. For a random vector $\vz\in
\mathbb{R}^p$, write $\vz\in{\cal L}^a$, $a > 0$, if $\| \vz
\|_a =: [ \E(|\vz|^a) ]^{1/a} < \infty$.

\section{Covariance matrix estimation for high-dimensional
stationary processes} \label{sec:stationary}

In this section we shall assume that $(\vz_i)$ is a
$p$-dimensional stationary process of the form
%
\begin{equation}
\label{eq:casual} \vz_i = \vg(\calF_i),
\end{equation}
where $\vg(\calF_i) = (g_1(\calF_i), \ldots, g_p(\calF_i))^\top$
is an $\mathbb{R}^p$-valued measurable function, $\calF_i =
(\ldots, \ve_{i-1}, \ve_i)$ is a shift process and $\ve_i$ are
i.i.d. random vectors. Following \mbox{\citet{MR2172215}}, we can view
$\calF_i$ and $\vz_i$ as the input and the output of a physical
system, respectively, and $\vg(\cdot)$ is the transform
representing the underlying physical mechanism. The framework
(\ref{eq:casual}) is quite general. Some examples are presented in
\citet{MR2812816}. It can also be conveniently extended to locally
stationary processes; see Section~\ref{sec:covariancematrixestiamtion_nonstatproc}.

Write $\vz_i = (Z_{1i}, \ldots, Z_{pi})^\top$ and $Z_{p \times n}
= (\vz_i)_{i=1}^n$, the data matrix observed at time points $i =
1,\ldots,n$. Here we shall consider estimation of the $p \times p$
covariance matrix $\Sigma= \operatorname{ cov}(\vz_i)$ based on the
realization $\vz_1, \ldots,\vz_n$, while Section~\ref{sec:precision_statproc} concerns estimation of its inverse.
We consider Frobenius and spectral norm convergence of the
\textit{thresholded estimator}
%
\begin{equation}
\label{eqn:thresholdestimator_covmat} T_u(\hat\Sigma_{n}) =\bigl(\hat{
\sigma}_{jk} \mathbb{I}\bigl(|\hat{\sigma}_{jk}| \ge u\bigr)
\bigr)_{1 \leq j,k \leq p},
\end{equation}
where $\hat\Sigma_n = (\hat{\sigma}_{jk})$ is the sample
covariance matrix defined in (\ref{eq:samplecovariancematrix});
see \citet{bickellevina2008a}. It was shown in the latter paper
that, with a properly chosen $u$, $T_u(\hat\Sigma_n)$ is a
consistent estimator when $\Sigma_0 \in\mathcal{G}_r(\tilde M)$
[see (\ref{eqn:strongell_qball})] and $(\vz_i)$ are i.i.d.
sub-Gaussian. Our rates of convergence depend on the dependence of
the process and the moment conditions, which can be quite mild.
Our main theoretical result is given in Section~\ref{subsec:covmat_stat_nonlinproc}. To obtain a consistent
estimate for $\Sigma$, we need to impose regularization
conditions. In particular, we shall assume that $\Sigma$ is weakly
dependent in that most of its entries are small, by providing a
bound on the tail empirical process of covariances. Some examples
are provided in Section~\ref{sec:stprocess} with applications to
spatial--temporal processes.

\subsection{Asymptotic results}
\label{subsec:covmat_stat_nonlinproc}

To establish a convergence theory for covariance matrix estimates,
we shall use the functional dependence measure of
\citet{MR2172215}. Recall that $Z_{ji} = g_j(\calF_i)$, $1 \le j
\le p$, where $g_j(\cdot)$ is the $j$th coordinate projection of
the $\mathbb{R}^p$-valued measurable function $\vg$. For $w > 0$,
the functional dependence measure of $Z_{ji}$ is defined by
%
\begin{equation}
\label{eq:functiondependencemeasure_stat} \theta_{i,w,j} =\bigl \|Z_{ji}-Z'_{ji}
\bigr\|_w = \bigl(\E\bigl|Z_{ji}-Z'_{ji}\bigr|^w
\bigr)^{1/w},
\end{equation}
where $Z'_{ji} = g_j(\calF_i')$, $\calF_i' = (\ldots, \ve_{-1},
\ve_0', \ve_1, \ldots, \ve_i)$ and $\ve_0'$ is such that $\ve_0',
\ve_l$, $l \in{\mathbb Z}$, are i.i.d. In other words, $Z'_{ji}$
is a coupled version of $Z_{ji}$ with $\ve_0$ in the latter
replaced by an i.i.d. copy $\ve_0'$. In \citet{MR2812816}
functional dependence measures were computed for some commonly
used linear and nonlinear stationary processes. We shall assume
that the short-range dependence (SRD) condition holds,
%
\begin{equation}
\label{eq:srdtail} \Theta_{m, w} = \max_{1\le j \le p} \sum
_{l=m}^\infty\theta_{l,w,j} < \infty.
\end{equation}
If (\ref{eq:srdtail}) fails, the process $(Z_{j i})_{i \in
{\mathbb Z}}$ may exhibit long-range dependence, and the asymptotic
behavior can be quite different. A nonlinear process satisfying
(\ref{eq:srdtail}) is given in Example~\ref{exmp:nonlin_stat},
while Example~\ref{exmp:linstat} concerns linear processes.
Theorems \ref{thmm:F08122} and \ref{thmm:spectral} provide rates of
convergence under the normalized Frobenius norm and the spectral
norm for the thresholded estimate $T_u(\hat{\Sigma}_n)$,
respectively. The constants~$C$ therein are independent of $n$,
$u$ and $p$.

\begin{thmm}
\label{thmm:F08122} Assume that there exist $q > 2$, $\alpha> 0$,
$\mu< \infty$ and a positive constant $C_0 < \infty$ such that
$\max_{j \le p} \|Z_{ji} \|_{2 q} \le\mu$ and $\Theta_{m, 2q} \le
C_0 m^{-\alpha}$ for all $m \ge1$. Let $\tilde\alpha= \alpha
\wedge(1/2-1/q)$ and $\tilde\beta= (3 + 2 \tilde\alpha q) /
(1+q)$. Define
%
\begin{eqnarray}
\label{eqn:D_u} H(u) &=& \cases{ %
u^{2-q}
n^{1-q}, & \quad$\mbox{if } \alpha> 1/2-1/q;$
\vspace*{2pt}\cr
u^{2-q} n^{1-q} (\log n)^{1+q}, & \quad$\mbox{if } \alpha=
1/2-1/q;$
\vspace*{2pt}\cr
u^{2-q} n^{-q(\alpha+1/2)}, & \quad$\mbox{if } \alpha< 1/2-1/q,$}\\
\label{eqn:D_u1} G(u) &=& \cases{ %
\bigl(n^{-1} + u^2\bigr) e^{-n u^2}, &\quad $\mbox{if }
\alpha> 1/2-1/q;$
\vspace*{2pt}\cr
\bigl(n^{-1} (\log n)^2 + u^2\bigr)
e^{-n (\log n)^{-2} u^2 }, &\quad $\mbox{if } \alpha= 1/2-1/q;$
\vspace*{2pt}\cr
\bigl(n^{-\tilde\beta} + u^2\bigr) e^{-n^{\tilde\beta} u^2}, &\quad $\mbox{if }
\alpha< 1/2-1/q$}
\end{eqnarray}
and
%
\begin{equation}
\label{eqn:F231038} D(u) = {1\over p^2} \sum_{j,k=1}^p
\bigl(u^2 \wedge\sigma_{jk}^2\bigr).
\end{equation}
Then there exists a constant $C$, independent of $u$, $n$ and $p$,
such that
%
\begin{equation}
\label{eq:A807148} {{\E|T_u(\hat{\Sigma}_n)-\Sigma|_F^2} \over p^2} \lesssim D(u) + \min \biggl(
{1\over n}, {u^{2-q} \over n^{q/2}}, H(u) + G(C u) \biggr).
\end{equation}
\end{thmm}

\begin{rmk}\label{rem:A031046}
If $\alpha> 1/2 - 1/q$, elementary calculations indicate that\break
$H(u) + G(C u) \lesssim u^{2-q} n^{-q/2}$. Hence the right-hand
side of (\ref{eq:A807148}) is $\asymp D(u) + \min(n^{-1},  H(u)
+ G(C u))$. The term $u^{2-q} n^{-q/2}$ is needed if $\alpha\le
1/2 - 1/q$. 
\end{rmk}

By Theorem~\ref{thmm:F08122}, if $u = O(n^{-1/2})$, then $p^{-2} \E
|T_u(\hat{\Sigma}_n)-\Sigma|_F^2 = O(n^{-1})$. Better convergence
rates can be achieved if $D(n^{-1/2}) = o(n^{-1})$ by choosing a
larger threshold; see cases (i)--(iii) in Corollary~\ref{cor:M140825} below.

\begin{cor}
\label{cor:M140825} Assume that the conditions of Theorem~\ref{thmm:F08122} hold. Let $\Upsilon= \inf_{u > 0} p^{-2} \E
|T_u(\hat{\Sigma}_n)-\Sigma|_F^2$; let $\tilde G(u) = \min(G(u),
u^{2-q} n^{-q/2})$ if $\alpha\le1/2-1/q$ and $\tilde G(u) =
G(u)$ if $\alpha> 1/2-1/q$.

Let $u_\diamond\ge n^{-1/2}$ be the unique solution to the
equation $H(u) = G(u)$. \textup{(i)} If $\bar{D} =: p^{-2} \sum_{j,k=1}^p
\sigma_{jk}^2 = O(H(1))$, then there is a fixed constant $c
> 0$ such that $\Upsilon\lesssim H(u) \asymp H(1)$ for all $u \in
[c, \mu]$. \textup{(ii)} If $H(1) = o(\bar{D})$ and $D(u_\diamond) \le
H(u_\diamond)$, let $u_\dagger$ solve $D(u_\dagger) =
H(u_\dagger)$, then $\Upsilon\lesssim D(u_\dagger)$. \textup{(iii)} If
$H(1) = o(\bar{D})$, $D(u_\diamond) > H(u_\diamond)$ and
$D(n^{-1/2}) =o(n^{-1})$, let $u_\circ$ be the solution to the
equation $D(u) = \tilde G(u)$ over the interval $u \in[n^{-1/2},
u_\diamond]$, then $\Upsilon\lesssim D(u_\circ)$. \textup{(iv)} If $n^{-1}
= O(D(n^{-1/2}))$, then the right-hand side of (\ref{eq:A807148})
is $\asymp n^{-1}$ for all $u \le n^{-1/2}$ and $\Upsilon\lesssim
n^{-1}$.
\end{cor}

Theorem~\ref{thmm:F08122} and Corollary~\ref{cor:M140825} describe
how the Frobenius rate of convergence depends on the sample size
$n$, the dimension $p$, the smallness measure quantified by the
function $D(u)$ and the heaviness of tails (moment conditions)
and strength of dependence which are characterized by $q$ and
$\alpha$, respectively. It suggests the interesting dichotomy
phenomenon: under the weaker dependence condition $\alpha>
1/2-1/q$, the thresholded estimate $T_u(\hat{\Sigma}_n)$ has the same
convergence rates as those obtained under independence. However,
the convergence becomes slower under stronger temporal dependence
with $\alpha< 1/2-1/q$. The phase transition occurring at $\alpha
= 1/2-1/q$. The theorem also provides information about the
optimal threshold $u$, as revealed in its proof. The optimal
threshold balances the bias or the smallness function $D(u)$, the
tail function $H(u)$ and the variance component which roughly
corresponds to the Gaussian-type function $G(u)$. Under different
conditions, the optimal threshold assumes different forms; see
Corollaries~\ref{cor:F_stationary} and \ref{cor:Fexp}.\vadjust{\goodbreak}

\begin{pf*}{Proof of Theorem~\ref{thmm:F08122}}
We first assume $\alpha> 1/2 - 1/q$. Note that
%
\begin{eqnarray}
\label{eqn:decompose_Frisk} \E\bigl|T_u(\hat{\Sigma}_n)-
\Sigma\bigr|_F^2 &=& \sum_{j,k=1}^p
\E\bigl[\hat{\sigma}_{jk}\mathbb{I}\bigl(|\hat{\sigma}_{jk}|\ge
u\bigr)-\sigma_{jk}\bigr]^2
\nonumber
\\[-8pt]
\\[-8pt]
\nonumber
&\le& 2\sum
_{j,k=1}^p \E\bigl(W_{j k}^2
\bigr) + 2B(u/2),
\end{eqnarray}
where $W_{j k} = \hat{\sigma}_{jk} \mathbb{I}(|\hat{\sigma}_{jk}|
\ge u)-\sigma_{jk}\mathbb{I}(|\sigma_{jk}|\ge u/2)$ and
%
\begin{equation}
\label{eqn:bias_part_general} B(u)=\sum_{j,k=1}^p
\sigma_{jk}^2\mathbb{I}\bigl(|\sigma_{jk}|< u\bigr).
\end{equation}
Let events $A^1_{jk} = \{|\hat{\sigma}_{jk}|\ge u, |\sigma_{jk}|
\ge u/2 \}$, $A^2_{jk} = \{|\hat{\sigma}_{jk}| < u, |\sigma_{jk}|
\ge u/2\}$ and $A^3_{jk} = \{|\hat{\sigma}_{jk}| \ge u,
|\sigma_{jk}| < u/2\}$, $1\le j, k \le p$. Observe that
\[
W_{j k} = W_{j k} \mathbb{I}\bigl(A^1_{jk}
\bigr) + W_{j k} \mathbb{I}\bigl(A^2_{jk}\bigr) +
W_{j k} \mathbb{I}\bigl(A^3_{jk}\bigr).
\]
We shall consider these three terms separately. Write $\xi_{j k} =
\hat{\sigma}_{jk} - \sigma_{jk}$.

\emph{Case} I: on the event $A^1_{jk}$,\vspace*{1pt} since the functional
dependence measure for the product process $Z_{j i} Z_{k i}$, $i
\in\mathbb Z$, satisfies
%
\begin{eqnarray}
\label{eq:Julypdm} \bigl\| Z_{j i} Z_{k i} - Z'_{j i}
Z'_{k i} \bigr\|_q &\le& \bigl\| Z_{j i}
Z_{k i} - Z'_{j i} Z_{k i}
\bigr\|_q + \bigl\| Z'_{j i} Z_{k i} -
Z'_{j i} Z'_{k i}
\bigr\|_q
\nonumber
\\[-8pt]
\\[-8pt]
\nonumber
&\le& \mu(\theta_{i, 2 q, j} + \theta_{i, 2 q, k}),
\end{eqnarray}
it follows from the moment inequality Theorem~2.1 in
\citet{MR2353389} that
%
\begin{equation}
\label{eq:F241256} \|\xi_{jk}\|_q \le c_q
n^{-1/2} \mu\Theta_{0, 2q},
\end{equation}
where $c_q$ is a constant only depending on $q$. Let $C_1 = c^2_q
\mu^2 \Theta_{0, 2q}^2$. Then
%
\begin{equation}
\label{eq:July131} \E\bigl\{ W_{j k}^2\mathbb{I}
\bigl(A^1_{jk}\bigr) \bigr\} \le\E \xi_{j k}^2
\mathbb{I}\bigl(|\sigma_{jk}|\ge u/2\bigr) \le C_1\frac{\mathbb{I}(|\sigma_{jk}|\ge u/2)}{n}.
\end{equation}

\emph{Case} II: on the event $A^2_{jk}$, we observe that
%
\begin{eqnarray}
\label{eq:July132} \E\bigl\{ W_{j k}^2\mathbb{I}
\bigl(A^2_{jk}\bigr) \bigr\} &= & \E\bigl[
\sigma_{jk}^2\mathbb{I}\bigl(|\sigma_{jk}|\ge u/2,|
\hat{\sigma }_{jk}| < u\bigr)\bigr]
\nonumber\\
&\le& 2\E\bigl[\xi_{jk}^2
\mathbb{I}\bigl(|\sigma_{jk}|\ge u/2,|\hat{\sigma}_{jk}| < u\bigr)
\bigr]
\nonumber
\\[-8pt]
\\[-8pt]
\nonumber
&&{} + 2 \E\bigl[\hat{\sigma}_{jk}^2 \mathbb{I}\bigl(|
\sigma_{jk}|\ge u/2,|\hat{\sigma}_{jk}| < u\bigr)\bigr]
\\
&\le&
2\bigl(C_1n^{-1}+u^2\bigr)\mathbb{I}\bigl(|
\sigma_{jk}|\ge u/2\bigr).\nonumber
\end{eqnarray}

\emph{Case} III: on the event $A^3_{jk}$, let
%
\begin{eqnarray}
\label{eqn:Delta_bound} \Delta_{jk} &=& {\E\bigl[\xi_{jk}^2
\mathbb{I}\bigl(|\hat{\sigma}_{jk}| \ge u, |\sigma_{jk}| < u/2\bigr)
\bigr]}
\nonumber\\
&=& \E\bigl[\xi_{jk}^2 \mathbb{I}\bigl(|\hat{
\sigma}_{jk}| \ge u, |\sigma_{jk}| < u/2, |
\xi_{jk}|>u/2\bigr)\bigr]
\\
&\le& \E\bigl[\xi_{jk}^2
\mathbb{I}\bigl(|\xi_{jk}|>u/2\bigr)\bigr].\nonumber
\end{eqnarray}
Then
%
\begin{eqnarray}
\label{eq:July647} \E\bigl\{ W_{j k}^2\mathbb{I}
\bigl(A^3_{jk}\bigr) \bigr\} &=& \E\bigl[\hat{
\sigma}_{jk}^2 \mathbb{I}\bigl(|\hat{\sigma}_{jk}|
\ge u, |\sigma_{jk}| < u/2\bigr)\bigr]
\nonumber
\\[-8pt]
\\[-8pt]
\nonumber
& \le& 2 \Delta_{jk} +
2\sigma_{jk}^2\mathbb{I}\bigl(|\sigma_{jk}| < u/2\bigr).
\end{eqnarray}
Since the functional dependence measure for the product process
$(Z_{j i} Z_{k i})_i$ satisfies (\ref{eq:Julypdm}), under the
decay condition $\Theta_{m, 2q} \le C m^{-\alpha}$, $\alpha>
1/2-1/q$, we have by Theorem~2(ii) in~\citet{liuxiaowu2012a} that
%
\begin{equation}
\label{eq:July14840} \Prob\bigl(|\xi_{j k}| > v\bigr) \le\frac{C_2n}{(n v)^q} +
C_3e^{-C_4n v^2}
\end{equation}
holds for all $v > 0$. Using integration by parts, we obtain
%
\begin{eqnarray}
\label{eq:A8051043} \E\bigl[\xi_{j k}^2 \mathbb{I}\bigl(|
\xi_{j k}|>v\bigr)\bigr] &=& v^2\Prob\bigl(|\xi_{j k}|>v\bigr) +
\int_{v^2}^\infty \Prob\bigl(|\xi_{j k}|>
\sqrt{w}\bigr) \,dw
\nonumber\\
&\le& v^2 \biggl[\frac{ C_2n}{(n v)^q} +
C_3e^{-C_4n v^2} \biggr]
\nonumber
\\[-8pt]
\\[-8pt]
\nonumber
& &{} + \int_{v^2}^\infty
\biggl[ \frac{C_2n}{(n\sqrt{w})^q} + C_3e^{-C_4nw} \biggr] \,dw
\\
&=&
C_5 n^{1-q} v^{2-q} + C_3
\bigl((C_4n)^{-1}+v^2\bigr)e^{-C_4n v^2},\nonumber
\end{eqnarray}
where $C_5 = C_2 q/(q-2)$. By (\ref{eq:F241256}), we also have
%
\begin{equation}
\label{eq:F24112} \E\bigl[\xi_{j k}^2 \mathbb{I}\bigl(|
\xi_{j k}|>v\bigr)\bigr] \le\min\biggl(\|\xi_{jk}
\|^2_2, {{\|\xi_{jk}\|^q_q} \over{v^{q-2} }}\biggr) \lesssim\min\biggl(
{1\over n}, {v^{2-q} \over n^{q/2}}\biggr).
\end{equation}

Combining cases I, II and III, by (\ref{eqn:bias_part_general}) and
(\ref{eq:July131})--(\ref{eq:F24112}), we have
%
\begin{eqnarray}
\label{eqn:Frisk_bound_general} {{\E|T_u(\hat{\Sigma}_n)-\Sigma|_F^2} \over p^2} &\lesssim& \frac{B(u/2)}{p^2} +
\frac{1 + nu^2}{np^2}\sum_{j,k=1}^p \mathbb{I}\bigl(|
\sigma_{jk}|\ge u/2\bigr)
\nonumber
\\[-8pt]
\\[-8pt]
\nonumber
& &{} + \min \biggl({1\over n},
{u^{2-q} \over n^{q/2}}, H(u) + G(C u) \biggr) =: M_0(u),
\end{eqnarray}
where $C = C_4^{1/2}/2$, and the constant of $\lesssim$ is
independent of $p$, $u$ and $n$. If $u \ge n^{-1/2}$, then
(\ref{eq:A807148}) clearly follows from the inequality $p^{-2}
\sum_{j,k} \mathbb{I}(|\sigma_{jk}|\ge v) \le v^{-2} D(v)$. If $u
< n^{-1/2}$, we also have (\ref{eq:A807148}) since in this case
$M_0(u) \asymp n^{-1}$ and the right-hand side of
(\ref{eq:A807148}) has the same order of magnitude $n^{-1}$.

The other cases with $0 < \alpha< 1/2-1/q$ and $\alpha= 1/2-1/q$
can be similarly handled. The key difference is that, instead of
(\ref{eq:July14840}), we shall now use the following versions of
Nagaev inequalities which can allow stronger dependence:
\begin{eqnarray*}
\label{eq:nagaevinequality_strongerdependence}
\Prob\bigl(|\xi_{jk}| > v\bigr) \le\cases{ %
\displaystyle \frac{C_2n^{q(1/2-\alpha)}}{(n v)^q} + C_3e^{-C_4n^{\tilde\beta} v^2},
 & \quad$\mbox{if } \alpha< 1/2-1/q;$
\vspace*{2pt}\cr
\displaystyle\frac{C_2 n(\log n)^{1+q} }{(n v)^q} + C_3 e^{-C_4n (\log n)^{-2} v^2},
&\quad
$\mbox{if } \alpha=
1/2-1/q.$}
\end{eqnarray*}
See also \citet{liuxiaowu2012a}.
\end{pf*}

\begin{pf*}{Proof of Corollary~\ref{cor:M140825}}
Let $M_1(u)$ be the term on the right-hand side of~(\ref{eq:A807148}). We now minimize $M_1(u)$ over $u > 0$. Let
\[
M_2(v) = D(v) + \min \biggl({1\over n},
{v^{2-q} \over n^{q/2}}, \max\bigl(H(v), G(v)\bigr) \biggr).
\]
Then $\inf_{u > 0} M_1(u) \asymp\inf_{v > 0} M_2(v)$. Clearly,
$\inf_{v \le n^{-1/2}} M_2(v) \asymp n^{-1}$. Let $v \ge
n^{-1/2}$. If $\alpha> 1/2-1/q$, then for some constant $c_q$, we
have $v^{2-q} n^{-q/2} \ge c_q v^2 e^{-n v^2} \ge c_q G(v) / 2$.
Also we have $v^{2-q} n^{-q/2} \ge H(v)$. Hence
%
\begin{equation}
\label{eqn:Fmax} \inf_{v \ge n^{-1/2}} M_2(v) \asymp\inf
_{v \ge n^{-1/2}} \max\bigl[D(v), H(v), G(v)\bigr].
\end{equation}
Note that the equation $H(u) = G(u)$ has a unique solution
$u_\diamond$ on $(n^{-1/2}, \infty)$, and the function $\max[H(u),
G(u)]$ is decreasing over $u \ge n^{-1/2}$. A plot of the
function in (\ref{eqn:Fmax}) is given in Figure~\ref{fig:M15608}(a). Let $u_\natural$ be the minimizer of the
right-hand side of~(\ref{eqn:Fmax}). For (i), assume $\bar{D} \le
C_0 n^{1-q}$ for some $C_0 > 0$. Then $u_\natural$ satisfies $D(u)
= H(u)$, which implies $u \ge C_0^{1/(2-q)}$, and hence (i)
follows. Note that (ii) follows in view of $u_\dagger= u_\natural
\ge u_\diamond$ and $u_\dagger\to0$. Similarly we have (iii)
since $u_\natural= u_\circ$. The last case (iv) is
straightforward since $M_2(u) \asymp n^{-1}$ for all $u \le
n^{-1/2}$.

If $0 < \alpha\le1/2-1/q$, assume $v \ge n^{-1/2}$, and then
(\ref{eqn:Fmax}) still holds with $G(v)$ therein replaced by
$\tilde G(v)$. A plot for this case is given in Figure~\ref{fig:M15608}(b). Note that $\tilde G(v) = G(v)$ if $v \ge
u_\diamond$. Then we can similarly have (i)--(iv).
\end{pf*}

\begin{rmk}
From the proof of Corollary~\ref{cor:M140825}, if $0 < \alpha\le
1/2-1/q$, in case~(iii), we can actually have the following
dichotomy: let $u_\triangle$ be the solution to the equation $G(u)
= u^{2-q} n^{-q/2}$. Then the minimizer $u_\natural\in[n^{-1/2},
u_\triangle]$ if $D(u_\triangle) \ge\tilde G(u_\triangle)$ and
$u_\natural\in[u_\triangle, u_\diamond]$ if $D(u_\triangle) \le
\tilde G(u_\triangle)$. For $\alpha> 1/2-1/q$, (\ref{eqn:Fmax})~indicates that $v^{2-q} n^{-q/2}$ is not needed; see also Remark~\ref{rem:A031046}. 
\end{rmk}

Using the argument for Theorem~\ref{thmm:F08122}, we can similarly
establish a spectral norm convergence rate.
\citet{bickellevina2008a} considered the special setting with
i.i.d. vectors. Our Theorem~\ref{thmm:spectral} is a significant
improvement by relaxing the independence assumption, by obtaining
a sharper rate and by presenting a moment bound. As in Theorem~\ref{thmm:F08122}, we also have the phase transition at $\alpha=
1/2 - 1/q$. Note that \citet{bickellevina2008a} only provides a
probabilistic bound.

\begin{thmm}
\label{thmm:spectral} Let the moment and the dependence conditions
in Theorem~\ref{thmm:F08122} be satisfied. Let $L_\alpha=
n^{1/q-1}, n^{1/q-1} (\log{n})^{1+1/q}, n^{-\alpha-1/2}$ and
$J_\alpha= n^{-1/2}, n^{-1/2} \log{n}, n^{-\tilde\beta/ 2}$, for
$\alpha> 1/2-1/q, \alpha= 1/2-1/q$, and $\alpha< 1/2-1/q$,
respectively. Define
%
\begin{equation}
D_\ast(u) = \max_{1 \le k \le p} \sum
_{j=1}^p \bigl(|\sigma_{jk}| \wedge u\bigr),\qquad
N_\ast(u) = \max_{1 \le k \le p} \sum
_{j=1}^p \I\bigl(|\sigma_{jk}| \ge u\bigr),
\end{equation}
and $ M_\ast(u) = L_\alpha p^{1/q} N^{1 + 1/q}_\ast(u) + J_\alpha
(\log{p})^{1/2} N_\ast(u)$. Then there exists a constant~$C$,
independent of $u$, $n$ and $p$, such that
%
\begin{eqnarray}\label{eq:F06626}
\bigl\|\rho\bigl(T_u(\hat\Sigma_n) - \Sigma\bigr)
\bigr\|_2 & \lesssim& D_\ast(u) + M_\ast(u/2)
\nonumber
\\[-8pt]
\\[-8pt]
\nonumber
 & &{} + p \min \biggl( {1 \over\sqrt{n}},
{u^{1-q/2} \over n^{q/4}}, \bigl(H(u) + G(Cu)\bigr)^{1/2} \biggr),
\end{eqnarray}
where $H(\cdot)$ and $G(\cdot)$ are given in (\ref{eqn:D_u}) and
(\ref{eqn:D_u1}), respectively.
\end{thmm}

\begin{pf}
We shall only deal with the weaker dependent case with $\alpha>
1/2-1/q$. The other cases similarly follow. Recall the proof of
Theorem~\ref{thmm:F08122} for $W_{jk}$, $\xi_{j k}$ and $A^l_{j k},
l =1, 2, 3$. Let matrices $V_l = (W_{jk} \I(A^l_{jk}))_{j, k \le
p}$. Similar to (\ref{eqn:bias_part_general}), let $B_\ast(u) =
\max_{1\le k \le p} \sum_{j=1}^p |\sigma_{jk}|
\I(|\sigma_{jk}|<u)$. Then
%
\begin{equation}
\label{eqn:F06524} \bigl|\rho\bigl(T_u(\hat\Sigma_n) - \Sigma
\bigr)\bigr| \le B_\ast(u/2) + \sum_{l=1}^3
\bigl|\rho(V_l)\bigr|.
\end{equation}
Let $N_k(u) = \{j\dvtx | \sigma_{j k}| \ge u/2 \}$ and $z_u = C_1
M_\ast(u/2)$, where $C_1 > 0$ is a large constant. Since
$\rho(V_1) \le\max_{k \le p} \sum_{j \in N_k(u)} |\hat\sigma_{jk}
- \sigma_{jk}| =: Q$, by (\ref{eq:July14840}),
%
\begin{eqnarray}
\label{eq:F02715p} { {\|\rho(V_1)\|^2_2} \over2} &\le& \int_0^\infty
z \Prob(Q \ge z) \,d z
\nonumber\\
&\lesssim& { {z_u^2} \over2} + \int
_{z_u}^\infty z p S_u \biggl[
{n \over{(n z/S_u)^q}} + e^{-C_4 n z^2 S_u^{-2}} \biggr] \,d z \\
&\lesssim&
M^2_\ast(u/2),\nonumber
\end{eqnarray}
where $S_u = N_\ast(u/2)$. Similar to (\ref{eq:July132}), since
$\sigma_{jk} \le|\hat\sigma_{jk} - \sigma_{jk}| + u$ on
$A^2_{jk}$,
%
\begin{equation}
\label{eqn:spectral_A2} \bigl|\rho(V_2)\bigr| \le Q + u S_u \le Q + 2
D_\ast(u).
\end{equation}
Using the idea of (\ref{eq:July647}), we have
%
\begin{eqnarray}
\label{eqn:spectral_A3} \rho^2(V_3) &\le&\sum
_{j, k}\bigl |W_{jk} \I\bigl(A^3_{jk}
\bigr)\bigr|^2
\nonumber
\\[-8pt]
\\[-8pt]
\nonumber
&\le&2 \sum_{j,k}
\xi_{jk}^2 \I\bigl(|\xi_{jk}| > u/2\bigr) + 2
B_\ast^2(u/2).
\end{eqnarray}
By (\ref{eqn:Delta_bound})--(\ref{eq:F24112}) and
(\ref{eqn:F06524})--(\ref{eqn:spectral_A3}), we have
(\ref{eq:F06626}) since $B_\ast(u/2) \le B_\ast(u) \le D_\ast(u)$.
\end{pf}

The bounds in Theorems \ref{thmm:F08122} and \ref{thmm:spectral}
depend on the smallness measures, the moment order $q$, the
dependence parameter $\alpha$, the dimension $p$ and the sample
size~$n$. The problem of selecting optimal thresholds is highly
nontrivial. Our numeric experiments show that the cross-validation
based method has a reasonably good performance. However, we are
unable to provide a theoretical justification of the latter
method, and pose it as an open problem.

\begin{exmp}[(Stationary Markov chains)]
\label{exmp:nonlin_stat} We consider
the nonlinear process~$(\vz_i)$ defined by the iterated random
function
%
\begin{equation}
\label{eq:iteratedrandomfun_stat} \vz_i=\vg(\vz_{i-1},\ve_i),
\end{equation}
where $\ve_i$'s are i.i.d. innovations, and $\vg(\cdot,\cdot)$ is
an $\mathbb{R}^p$-valued and jointly measurable function, which
satisfies the following two conditions: (i) there exists some
$\vx_0$ such that $\|\vg(\vx_0, \ve_0)\|_{2q}<\infty$ and (ii)
%
\begin{equation}
\label{eqn:stochatic_Lip_const_stationary}
L = \sup_{\vx\neq\vx'} \frac{\|\vg(\vx,\ve_0)-\vg(\vx',\ve_0)
\|_{2q}}{|\vx-\vx'|} < 1.
\end{equation}
Then, it can be shown that $\vz_i$ defined in
(\ref{eq:iteratedrandomfun_stat}) has a stationary ergodic
distribution $\vz_0 \in{\cal L}^{2q}$ and, in addition, $(\vz_i)$
has the \emph{geometric moment contraction} (GMC) property; see
\citet{wushao2004a} for details. Therefore, we have $\Theta_{m,2q}
= O(L^m)$ and Theorems \ref{thmm:F08122} and \ref{thmm:spectral}
with $\alpha> 1/2-1/q$ and $\tilde\beta=1$ can be applied.
\end{exmp}

\begin{exmp}[(Stationary linear processes)]
\label{exmp:linstat} An important
special class of (\ref{eq:casual}) is the vector linear process
%
\begin{equation}
\label{eq:linearprocess} \vz_i = \sum_{m=0}^\infty
A_m \ve_{i-m},
\end{equation}
where $A_m, m \ge0$, are $p \times p$ matrices, and $\ve_i$ are
i.i.d. mean zero random vectors with finite covariance matrix
$\Sigma_\ve= \E(\ve_i \ve_i^\top)$. Then $\vz_i$ exists almost
surely with covariance matrix $\Sigma= \sum_{m = 0}^\infty A_m
\Sigma_\ve A_m^\top$ if the latter converges. Assume that the
innovation vector $\ve_i = (e_{1 i}, \ldots, e_{p i} )^\top$,
where $e_{ji}$ are i.i.d. with mean zero, variance $1$ and $e_{ji}
\in{\cal L}^{2q}$, $q > 2$, and the coefficient matrices $A_i =
(a_{i, j k})_{1\le j, k \le p}$ satisfy $\max_{j \le p}
\sum_{k=1}^p a_{i, j k}^2 = O(i^{-2-2\gamma})$, $\gamma> 0$. By\vadjust{\goodbreak}
Rosenthal's inequality, the functional dependence measure
$\theta^2_{i, 2 q, j} \le c_q \sum_{k=1}^p a_{i, j k}^2 =
O(i^{-2-2\gamma})$, and hence by (\ref{eq:srdtail}) $\Theta_{m, 2
q} = O(m^{-\gamma})$. By Theorem~\ref{thmm:F08122}, the normalized
Frobenius norm of the thresholded estimator has a convergence rate
established in (\ref{eq:A807148}) with $\alpha= \gamma$, $\tilde
\alpha= \gamma\wedge(1/2-1/q)$ and $\tilde\beta$. Note that
our moment condition relaxes the commonly assumed sub-Gaussian
condition in previous
literature~[\citet{rothmanbickellevinazhu2008a,lamfan2009a,zhoulaffertywasserman2010a}].
For the vector AR(1) process $\vz_i =
A \vz_{i-1} + \ve_i$, where $A$ is a real matrix with spectral
norm $\rho(A) < 1$, it is of form (\ref{eq:linearprocess}) with
$A_m = A^m$, and the functional dependence measure $\theta_{i, 2
q, j} = O(\rho(A)^i)$. The rates of convergence established in
(\ref{eq:A807148}) hold with $\alpha> 1/2-1/q$ and $\tilde\beta
= 1$.
\end{exmp}

\subsection{Positive-definitization}
The thresholded estimate $T_u(\hat\Sigma_n)$ may not be positive
definite. Here we shall propose a simple modification that is
positive definite and has the same rate of convergence. Let
$T_u(\hat\Sigma_n) = \mathbf{Q} \hat\Lambda\mathbf{Q}^\top=
\sum_{j=1}^p \hat\lambda_j \vq_j\vq_j^\top$ be its
eigen-decomposition, where $\mathbf{Q}$ is an orthonormal matrix and
$\hat\Lambda$ is a diagonal matrix. For $v > 0$, consider
%
\begin{equation}
\label{eqn:modified_threshold_estimator} \tilde S_v = \sum_{j=1}^p
(\hat{\lambda}_j \vee v) \vq_j\vq_j^\top,
\end{equation}
where $0 < v \le\sqrt p \varpi$ and $\varpi^2$ is the rate of
convergence in (\ref{eq:A807148}). Let $\mu_1, \ldots, \mu_p$ be
the diagonal elements of $\mathbf{Q}^\top\Sigma\mathbf{Q}$. Then we
have by Theorem~\ref{thmm:F08122} that $\sum_{j=1}^p
(\hat{\lambda}_j - \mu_j)^2 \le p^2 \varpi^2$, and consequently
\begin{eqnarray*}
|\tilde S_v - \Sigma|_F^2 &\le& 2 \bigl|
\tilde S_v - T_u(\hat\Sigma_n)\bigr|_F^2
+ 2 \bigl|T_u(\hat\Sigma_n) - \Sigma\bigr|_F^2
\\
&\le& 2 \sum_{j=1}^p \bigl(\hat{
\lambda}_j - (\hat{\lambda}_j \vee v)
\bigr)^2 + 2 \varpi^2 p^2
\\
&\le& 2 \sum
_{j=1}^p \bigl(2 \hat{
\lambda}_j^2 \mathbf{1}_{\hat{\lambda}_j \le0} +
2v^2\bigr) + 2 \varpi^2 p^2.
\end{eqnarray*}
If $\hat{\lambda}_j \le0$, since $\mu_i \ge0$, we have
$|\hat{\lambda}_j| \le|\hat{\lambda}_j - \mu_i|$. Then $|\tilde
S_v - \Sigma|_F^2 \le4 v^2 p + 6 \varpi^2 p^2 \le10 \varpi^2
p^2$. Note that the eigenvalues of $\tilde S_v$ are bounded below
by $v$, and thus it is positive definite. In practice we suggest
using  $v =\break (p^{-1} \sum_{j, k=1}^p u^2\times  \I(|\hat\sigma_{j k}| \ge
u))^{1/2}$. The same positive-definization procedure also applies
to the spectral norm and its rate can be similarly preserved.

\subsection{Classes of covariance matrices}
\label{sec:stprocess} In this section we shall compute the
smallness measure $D(u)$ for certain class of covariance matrices,
so that Theorem~\ref{thmm:F08122} is applicable. We consider some
widely used spatial processes. Let the vectors $\vz_i = (Z_{1i},
\ldots, Z_{pi})^\top$, $1 \le i \le n$, be observed at sites
$s_1^\circ, \ldots, s_p^\circ\in{\mathbb R^2}$. Assume that the\vadjust{\goodbreak}
covariance function between $Z_{j i}$ and $Z_{k i}$ satisfies
%
\begin{equation}
\sigma_{j k} = \operatorname{ cov}(Z_{j i}, Z_{k i}) = f
\bigl( d\bigl(s_j^\circ, s_k^\circ
\bigr) \bigr),
\end{equation}
where $ d(s_j^\circ, s_k^\circ) $ is a distance between sites
$s_j^\circ$ and $s_k^\circ$, and $f$ is a real-valued function with
$f(0) = 1$ and $f(\infty) = 0$. For example, we can choose $d(s,
s') = |s - s'|$ as the Euclidean distance between sites $s$ and
$s'$. Assume that, as $m \to\infty$,
%
\begin{equation}
\label{eq:A811032} f(m) = O\bigl(m^{-K}\bigr),
\end{equation}
where the index $K > 0$ characterizes the spatial dependence, or
%
\begin{equation}
\label{eq:A821010} f(m) \le\exp\bigl(-C (m/\tau)^\theta\bigr), \qquad 0 < \theta
\le2,
\end{equation}
where $\tau$ is the characteristic length-scale, and
%
\begin{equation}
\label{eq:A811028} {1\over p} \sum_{j, k = 1}^p
\mathbb{I}\bigl(d\bigl(s_j^\circ, s_k^\circ
\bigr) \le m\bigr) = O\bigl(m^\chi\bigr).
\end{equation}
Condition (\ref{eq:A811028}) outlines the geometry of the sites
$(s_j^\circ)_{j=1}^p$, and $\chi$ can be roughly interpreted as
the correlation dimension. It holds with $\chi= 2$ if $s_j^\circ$
are ${\mathbb Z^2}$ points in a disk or a square, and $\chi= 1$
if $s_j^\circ= (j, 0), j=1, \ldots, p$. The rational quadratic
covariance function [\citet{MR2514435}] is an example of
(\ref{eq:A811032}), and it is widely used in spatial statistics,
%
\begin{equation}
\label{eqn:rationalquad_covfuns} f(m) = \biggl(1+\frac{m^2}{K \tau^2} \biggr)^{-K/2},
\end{equation}
where $K$ is the smoothness parameter and $\tau>0$ is the length
scale parameter. We now provide a bound for $D(u)$. By
(\ref{eq:A811032}) and (\ref{eq:A811028}), as $u \downarrow0$,
the covariance tail empirical process function
%
\begin{equation}
\label{eqn:A02736} F(u) =: {1\over p^2} \sum
_{j,k=1}^p \mathbb{I}\bigl(|\sigma_{jk}| \ge u\bigr)
\le p^{-1} \min\bigl(p, C u^{-\chi/K}\bigr)
\end{equation}
for some constant $C > 0$ independent of $n$, $u$ and $p$. If $K >
\chi/2$, then
%
\begin{eqnarray}
\label{eqn:A02849} D(u) &=& u^2 F(u) + \frac{1}{p^2} \sum
_{l=0}^\infty \sum_{j, k = 1}^p
\sigma_{jk}^2 \mathbb{I}\bigl( u 2^{-l-1} \le|
\sigma_{jk}| < u 2^{-l}\bigr)
\nonumber\\
&\le& u^2 F(u) +
\sum_{l=0}^\infty \bigl(2^{-l} u
\bigr)^2 F\bigl(2^{-l-1} u\bigr)
\\
&\le& u^2
p^{-1} \min\bigl(p, C u^{-\chi/K}\bigr) = u^2 \min
\bigl(1, C p^{-1} u^{-\chi/K}\bigr).\nonumber
\end{eqnarray}
In the strong spatial dependence case with $K < \chi/2$, we have
%
\begin{equation}
D(u) \le\min\bigl(C p^{-K}, u^2\bigr).
\end{equation}
To this end, it suffices to prove this relation with $u^2 >
p^{-K}$. Let $u_0 = p^{-K/\chi}$. Then
\begin{eqnarray*}
\bar{D} &\le& \sum_{l=0}^\infty
\bigl(2^{1+l} u_0\bigr)^2 F
\bigl(2^l u_0\bigr)
\cr
&\le& \sum
_{l=0}^\infty\bigl(2^{1+l} u_0
\bigr)^2 C p^{-1} \bigl(2^{1+l} u_0
\bigr)^{-\chi/K} \le C p^{-2K/\chi}.
\end{eqnarray*}

Class (\ref{eq:A821010}) allows the $\gamma$-exponential
covariance function with $f(m) =  \exp(-(m/\tau)^\gamma)$, and
some Mat\'ern covariance functions [\citet{MR1697409}] that are
widely used in spatial statistics. With (\ref{eq:A811028}),
following the argument in (\ref{eqn:A02849}), we can similarly
have
%
\begin{equation}
D(u) \le\min\bigl(u^2, C p^{-1} \tau^\chi
u^2\bigl(\log\bigl(2+u^{-1}\bigr)\bigr)^{\chi/\theta}
\bigr).
\end{equation}

Corollary~\ref{cor:F_stationary} of Theorem~\ref{thmm:F08122}
concerns covariance matrices satisfying (\ref{eqn:A02736}).
Slightly more generally, we introduce a decay condition on the
tail empirical process of covariances. Note that
(\ref{eqn:A02736}) is a special case of (\ref{eqn:sparsity_def})
with $M = C p$ and $r = \chi/ K$. For
(\ref{eqn:rationalquad_covfuns}) with possibly large length scale
parameter $\tau$, we can let $M = C \tau^2 p$. Similarly,
Corollary~\ref{cor:Fexp} can be applied to $f$ satisfying
(\ref{eq:A821010}) and the class ${\cal L}_r(M)$ defined in
(\ref{eqn:Expsparsity}), with $M = p \tau^\chi$ and $r = \chi/
\theta$.

\begin{defn}
\label{def:csc} For $M > 0$, let ${\cal H}_r(M)$, $0 \le r < 2$,
be the collection of $p \times p$ covariance matrices $\Sigma=
(\sigma_{jk})$ such that $\sup_{j \le p} \sigma_{jj} \le1$ and,
for all $0 < u \le1$,
%
\begin{equation}
\label{eqn:sparsity_def} \sum_{j,k=1}^p \mathbb{I}\bigl(|
\sigma_{jk}| \ge u\bigr) \le M u^{-r},
\end{equation}
and ${\cal L}_r(M)$, $r > 0$, be the collection of $\Sigma=
(\sigma_{jk})$ with $\sup_{j \le p} \sigma_{jj} \le1$ and
%
\begin{equation}
\label{eqn:Expsparsity} \sum_{j,k=1}^p \mathbb{I}\bigl(|
\sigma_{jk}| \ge u\bigr) \le M \log^r\bigl(2+u^{-1}
\bigr).
\end{equation}
\end{defn}

\begin{cor}
\label{cor:F_stationary} Assume (\ref{eqn:sparsity_def}). Let
conditions in Theorem~\ref{thmm:F08122} be satisfied and $\alpha>
1/2-1/q$. Let $\Upsilon= p^{-2} \sup_{\Sigma\in{\cal H}_r(M)}
\inf_{u > 0} \E|T_u(\hat{\Sigma}_n)-\Sigma|_F^2$. \textup{(i)} If $n^{q-1}
= O(p^2/M)$, then for $u \asymp1$, $\Upsilon= O(H(u)) =
O(n^{1-q})$. \textup{(ii)} If $p^2/M = o(n^{q-1})$ and $n^{(r+q)/2-1} (\log
n)^{(q-r)/2} \le p^2/M$, let $u'_\dagger= (n^{1-q} p^2/M
)^{1/(q-r)}$, then $\Upsilon= O({u'_\dagger}^{2-q} n^{1-q})$.
\textup{(iii)} If $p^2/M = o(n^{q-1})$ and
%
\begin{equation}
{n^{1-q/2} \over(\log n)^{(q-r)/2}} \le{M \over p^2} n^{r/2} \le1,
\end{equation}
then the equation $u^{2-r} M / p^2 = u^2 e^{-n u^2}$ has solution
$u'_\circ\asymp[n^{-1} \log(2 +\break   p^2 M^{-1} n^{-r/2})]^{1/2}$ and
$\Upsilon= O({u'_\circ}^{2-r} M/p^2)$. \textup{(iv)} If $n^{r/2} \ge p^2 /
M$, then the right-hand side of (\ref{eq:A807148}) is $\asymp
n^{-1}$ for $u=O(n^{-1/2})$ and $\Upsilon= O(n^{-1})$.\vadjust{\goodbreak}

In particular, if $p^2/M \asymp n^\phi$, $\phi> 0$, then we have
\textup{(i), (ii), (iii)} or \textup{(iv)} if $\phi> q-1$, $q-1 > \phi>
(q+r-2)/2$, $(q+r-2)/2 > \phi> r/2$ or $r/2 > \phi$ holds,
respectively.
\end{cor}

\begin{pf}
Similar to (\ref{eqn:A02849}), we have $D(u) \le\min(u^2, C
u^{2-r} M/p^2)$. Note that the solution $u_\diamond\ge n^{-1/2}$
to the equation $H(u) = G(u)$ satisfies $u_\diamond\sim( (q/2-1)
n^{-1} \log n)^{1/2}$. Then by Corollary~\ref{cor:M140825},
(i)--(iv) follow from elementary but tedious manipulations. Details
are omitted.
\end{pf}

By taking into consideration of $M$ in the tail empirical process
condition (\ref{eqn:sparsity_def}), we can view $p^2/M$ as the
\textit{effective dimension}. Corollary~\ref{cor:F_stationary}
describes the choice of the optimal threshold $u$ at different
regions of the effective dimension $p^2/M$ and the sample size
$n$. Case (i) [resp., (iv)] corresponds to the overly large (resp.,
small) dimension case. The most interesting cases are (ii) and
(iii). For the former, the tail function $H(\cdot)$ determines the
rate of convergence with a larger threshold~$u_\dagger$, while for
the latter with moderately large dimension the Gaussian-type
function $G(\cdot)$ leads to the optimal threshold $u_\circ<
u_\dagger$.

\begin{cor}
\label{cor:Fexp} Assume (\ref{eqn:Expsparsity}). Let conditions in
Theorem~\ref{thmm:F08122} be satisfied with $\alpha> 1/2-1/q$ and
$\Upsilon= p^{-2} \sup_{\Sigma\in{\cal L}_r(M)} \inf_{u > 0} \E
|T_u(\hat{\Sigma}_n)-\Sigma|_F^2$. \textup{(i)}~If $n^{q-1} = O(p^2/M)$,
then for $u \asymp1$, $\Upsilon= O(H(u)) = O(n^{1-q})$. \textup{(ii)} If
$p^2/M = o(n^{q-1})$ and $n^{q/2-1} (\log n)^{r+q/2} \le p^2/M$,
let $\varepsilon_\dagger= n^{1-q} p^2 / M$ and $u'_\dagger=\break
\varepsilon_\dagger^{1/q} (\log(2+\varepsilon_\dagger^{-1}))^{-r/q}$.
Then $\Upsilon= O({u'_\dagger}^{2-q} n^{1-q})$. \textup{(iii)} If
$n^{q/2-1} (\log n)^{r+q/2} > p^2/M \ge(\log n)^r$, let $\eta=
(\log n)^{-r} p^2/M$. If $\eta\ge2^{-r}$ let $u'_\circ= (n^{-1}
\log\eta)^{1/2}$. Then $\Upsilon= O(n^{-1} \eta^{-1} \log
\eta)$. \textup{(iv)} If $\eta$ in \textup{(iii)} is less than $2^{-r}$, then the
right-hand side of~(\ref{eq:A807148}) is $\asymp n^{-1}$ for
$u=O(n^{-1/2})$ and $\Upsilon= O(n^{-1})$.
\end{cor}

\begin{pf} We have $D(u) = u^2 \min(1, p^{-2} M
\log^r(2+u^{-1}))$. We shall again apply Corollary~\ref{cor:M140825}. Case (i) is straightforward. For (ii), we note
that the equation $u^q \log^r(2+u^{-1}) = \varepsilon$ has solution
$u_\dagger\asymp\varepsilon_\dagger^{1/q} (\log
(2+\varepsilon_\dagger^{-1}))^{-r/q}$. Under (iii), the equation $u^2
p^{-2} M \log^r(2+u^{-1}) = G(u)$ has solution $\asymp u'_\circ$.
\end{pf}

Corollaries \ref{cor:F_stationary} and \ref{cor:Fexp} deal with
the weaker dependence case with $\alpha> 1/2-1/q$. By Corollary~\ref{cor:M140825}, similar versions can be obtained for $\alpha
\le1/2-1/q$. Details are omitted.

\begin{figure}[t!]
\centering
\begin{tabular}{@{}c@{}}

\includegraphics{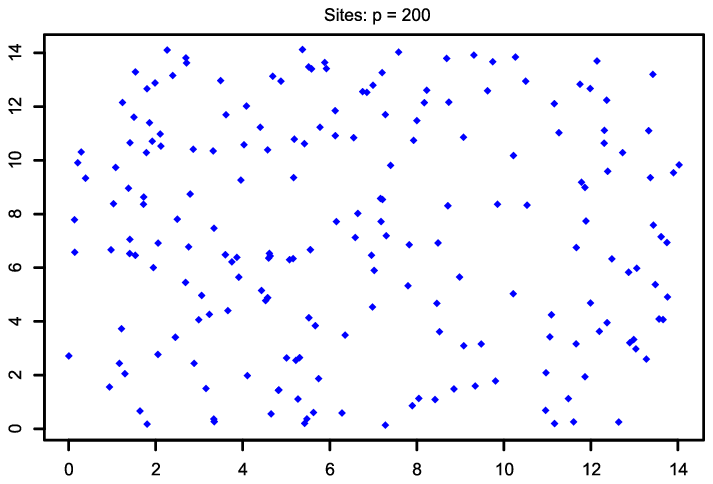}
 \\
\footnotesize{(a) $p$ sites $s_1^\circ, \ldots, s_p^\circ$
uniformly sampled from the square $p^{1/2} \times p^{1/2}$}
\end{tabular}\vspace*{4pt}
\begin{tabular}{@{}c@{\hspace*{2pt}}c@{}}

\includegraphics{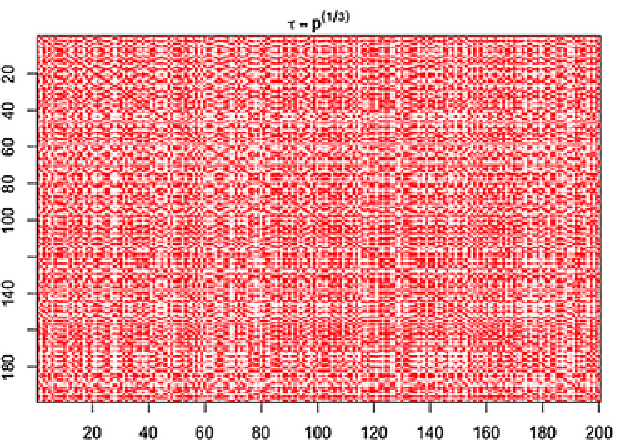}
&
\includegraphics{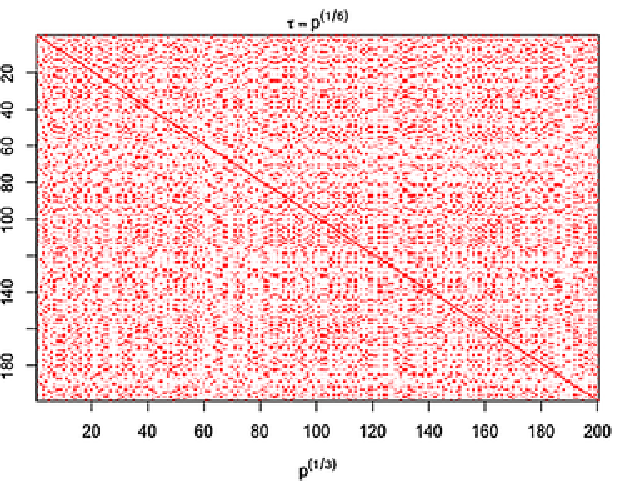}
\\
\footnotesize{(b) $\Sigma$: $\tau=p^{1/3}$}&
\footnotesize{(c) $\Sigma$:
$\tau=p^{1/6}$}
\end{tabular}\vspace*{4pt}
\begin{tabular}{@{}c@{}}

\includegraphics{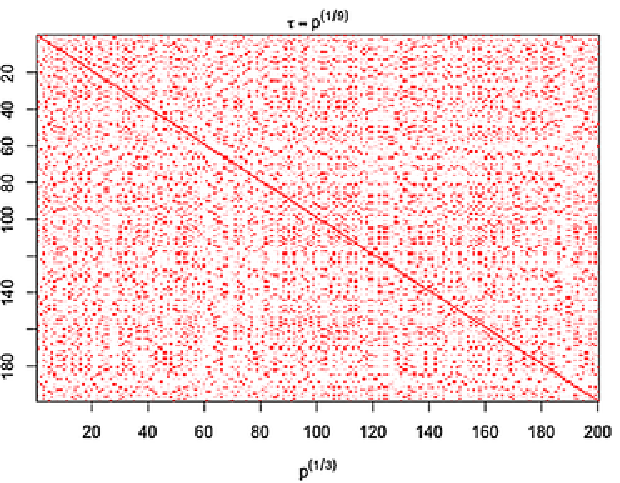}
\\
\footnotesize{(d) $\Sigma$:
$\tau=p^{1/9}$}
\end{tabular}
\caption{Rational quadratic covariance matrix $\Sigma$ for the
uniform random sites model on the $[0, p^{1/2}]^2$
square with three different scale length parameters:
$\tau= p^{1/3}, p^{1/6}$ and
$p^{1/9}$.}\label{fig:rational_quad_cov_mat}\vspace*{-4pt}
\end{figure}

\begin{figure}
\centering
\begin{tabular}{@{}c@{}}

\includegraphics{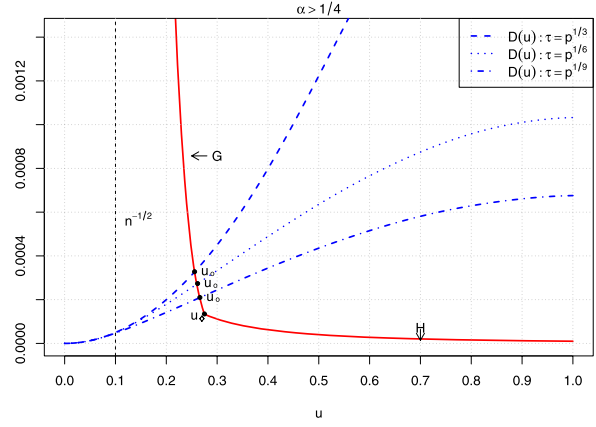}
\\
\footnotesize{(a) Weaker temporal dependence with
$\alpha> 1/4$}\\[3pt]

\includegraphics{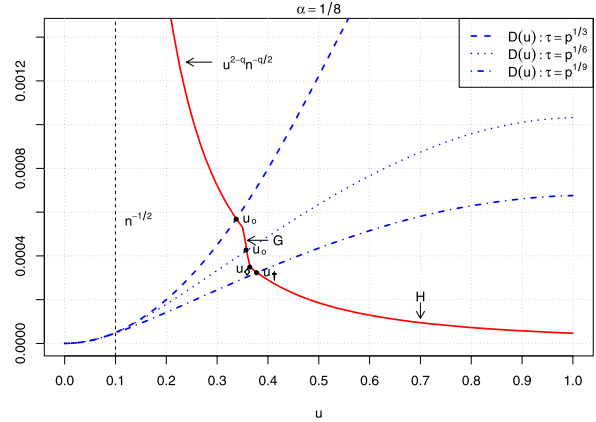}
\\
\footnotesize{(b) Stronger temporal dependence with
$\alpha= 1/8$}
\end{tabular}
\caption{Rates of convergence for the thresholded estimator in the
weaker ($\alpha> 1/4$) and stronger ($\alpha= 1/8$)
temporal dependence cases.}
\label{fig:M15608}
\end{figure}

As a numeric example, we use the rational quadratic covariances
(\ref{eqn:rationalquad_covfuns}) to illustrate the rates of
convergence given in Theorem~\ref{thmm:F08122} and Corollary~\ref{cor:M140825}. We choose $n = 100$, $p = 200$, $K = 4$, the
moment $q=4$ and consider the weaker ($\alpha> 1/4$) and
stronger ($\alpha= 1/8$) temporal dependence cases. We first
generate $p$ random sites uniformly distributed on the $p^{1/2}
\times p^{1/2}$ square; see
Figure~\ref{fig:rational_quad_cov_mat}(a).
Figure~\ref{fig:rational_quad_cov_mat}(b),
\ref{fig:rational_quad_cov_mat}(c) and
\ref{fig:rational_quad_cov_mat}(d) show three $200 \times200$
rational quadratic covariance matrices
(\ref{eqn:rationalquad_covfuns}) respectively with length scale
parameters $\tau= p^{1/3}, p^{1/6}$ and $p^{1/9}$, which
correspond to different levels of spatial dependence. Next, we
calculate the terms in Corollary~\ref{cor:M140825} for the
thresholded estimator. The results are shown in
Figure~\ref{fig:M15608}. In the plots, $u_\diamond$ is the
solution of $G(u)=H(u)$. Note that, $u_\natural$, the minimizer of
$\max[D(u), H(u), G(u)]$ over $u \ge n^{-1/2}$, can be
either $u_\dagger$ or $u_\circ$. We observe that when the spatial
dependence decreases, that is, the covariance matrix $\Sigma$ has more
small entries [e.g., Figure~\ref{fig:rational_quad_cov_mat}(d)], a
larger threshold is needed to yield the optimal rate of
convergence. When the temporal dependence increases (i.e., $\alpha
= 1/8$), a~larger threshold is needed and the rate of convergence
is slower than the one in the weaker dependence case (i.e., $\alpha
> 1/4$).

%
%

\subsection{Comparison with earlier results}

We now compare (\ref{eqn:sparsity_def}) with the commonly used
sparsity condition defined in terms of the \emph{strong
$\ell^q$-ball}
[\citet{bickellevina2008a,caizhou2011a,cailiuluo2011a}]
%
\begin{equation}
\label{eqn:strongell_qball} \mathcal{G}_r(\tilde{M}) = \Biggl\{\Sigma \Big|\max
_{j \le p} \sigma_{jj} \le1; \max_{1\le k \le
p}
\sum_{j=1}^p |\sigma_{jk}|^r
\le\tilde{M}\Biggr\},\qquad  0 \le r < 1.
\end{equation}
When $r = 0$, (\ref{eqn:strongell_qball}) becomes $\max_{1 \le k
\le p}\sum_{j=1}^p \mathbb{I}(\sigma_{jk}\neq0) \le\tilde{M}$, a
sparsity condition in the rigid sense. We observe that condition
(\ref{eqn:sparsity_def}) defines a broader class of sparse
covariance matrices in the sense that $\mathcal{G}_r(M/p)
\subset{\cal H}_r(M)$, which follows from
\[
\sum_{j,k} \I\bigl(|\sigma_{jk}|\ge u\bigr) \le p
\max_k\sum_j
\frac{|\sigma_{jk}|^{r}}{u^r} \le M u^{-r}.
\]
Hence Corollary~\ref{cor:F_stationary} generalizes the consistency
result of $T_u(\hat{\Sigma}_n)$ in \citet{bickellevina2008a} to the
non-Gaussian time series. Note that our convergence is in ${\cal
L}^2$ norm, while the error bounds in previous work [see, e.g.,
\citeauthor{bickellevina2008a}
(\citeyear{bickellevina2008a,bickellevina2008b})] are of probabilistic
nature; namely in the form $|T_u(\hat{\Sigma}_n) - \Sigma|_F^2$ is
bounded with large probability under the strong $\ell^q$-ball
conditions.

The reverse inclusion ${\cal H}_r(M)\subset\mathcal{G}_r(M/p)$ may
be false since the class $\mathcal{G}_r$ specifies the uniform
size of sums in matrix columns, whereas (\ref{eqn:sparsity_def})
can be viewed as an overall smallness measure over all entries of
the matrix. As an example, consider the covariance matrix
%
\begin{equation}
\label{eqn:counterexample} \Sigma_{p \times p} =
\pmatrix{ 1 & \varepsilon& \varepsilon&
\cdots& \varepsilon\vspace*{2pt}
\cr
\varepsilon& 1 & 0 & \cdots& 0 \vspace*{2pt}
\cr
\varepsilon& 0 & 1 & \cdots& 0 \vspace*{2pt}
\cr
\vdots& \vdots& \vdots&
\ddots& \vdots\vspace*{2pt}
\cr
\varepsilon& 0 & 0 & \cdots& 1 }
\end{equation}
where $0 < \varepsilon\le(p-1)^{-1/2}$ so that $\Sigma$ is
positive-definite. Then for any threshold level $u \in
(\varepsilon,1)$, $\sum_{j,k=1}^p \mathbb{I}(|\sigma_{jk}|\ge u) =
p$ and for any $u \in(0, \varepsilon]$,\break  $\sum_{j,k=1}^p
\mathbb{I}(|\sigma_{jk}|\ge u) = 3p-2$. In both cases, we may
choose $M=O(p)$. On the other hand, $\max_k \sum_j |\sigma_{jk}|^r
= 1+(p-1)\varepsilon^r$. So $\Sigma\notin\mathcal{G}_r(M/p)$ for
any $\varepsilon\ge(p-1)^{\eta/r-1/r}$ with $\eta\in(0,1-r/2)$.

With the strong $\ell^q$-ball and sub-Gaussian conditions,
\citet{caizhou2011a} showed that the minimax rate under the Bregman
divergence is $O(n^{-1} + \tilde{M}(\log{p}/n)^{1-r/2})$.
Observing that the upper bounds in Corollary~\ref{cor:F_stationary} is established under the larger parameter
space $ {\cal H}_r(M)\supset\mathcal{G}_r(\tilde{M})$ where
$M=p\tilde{M}$ and milder polynomial moments conditions, the lower
bound of \citet{caizhou2011a} automatically becomes a lower bound
in our setup. Therefore, in the moderately high-dimensional
situation with weaker temporal dependence, we can conclude that
the Frobenius norm bound in Corollary~\ref{cor:F_stationary}(iii)
is minimax rate optimal.

\begin{cor}
\label{cor:Frisk_bound_minimax_stationary_covmat} Let $\alpha
> 1/2-1/q$. Under the conditions in Corollary~\ref{cor:F_stationary}\textup{(iii)} and in addition assume $p^2 M^{-1}
n^{-r/2} \ge p^{\varepsilon}$ for some $\varepsilon>0$. Then
%
\begin{equation}
\label{eqn:Frisk_bound_minimax_stationary_covmat}
\inf_{\hat\Sigma} \sup_{\Sigma\in{\cal H}_r(M)}
p^{-1}\E |\hat{\Sigma}-\Sigma|_F^2 \asymp
\frac{M}{p} \biggl(\frac{\log{p}}{n} \biggr)^{1-{r}/{2}},
\end{equation}
where the inf is taken over all possible estimators based on the
data $Z_{p\times n}$.
\end{cor}

We next compare our Theorem~\ref{thmm:spectral} with the result in
Section~2.3 of \citet{bickellevina2008a}, where the special class
(\ref{eqn:strongell_qball}) is considered. Assuming $\max_{j}
\|Z_{ji}\|_{2q} \le\mu$, they obtained the \textit{probabilistic
bound}
%
\begin{equation}
\label{eq:F07934} \rho\bigl(T_{u_\BL}(\hat\Sigma_n) - \Sigma
\bigr) = O_p\bigl(\tilde M u_\BL^{1-r}\bigr)\qquad
\mbox{where } u_{\BL} = C p^{2/q} n^{-1/2},
\end{equation}
and $C > 0$ is a sufficiently large constant. As a natural
requirement for consistency, we assume $u_{\BL} \to0$, namely $p
= o(n^{q/4})$. Since $\Sigma\in{\cal G}_r(\tilde M)$, we have
$D_\ast(u) \le\tilde M u^{1-r} =: \bar D_\ast(u)$ and $N_\ast(u)
\le\min(p, \tilde M u^{-r}) =: \bar N_\ast(u)$. Consider the
weaker dependence case with $\alpha> 1/2-1/q$. Note that in
(\ref{eq:F06626}) $D_\ast(\cdot)$ is nondecreasing, while all
other three functions are nonincreasing. Let $u_1$, $u_2$, $u_3$
be the solutions to the equations $\bar N^{1+1/q}_\ast(u) p^{1/q}
n^{1/q-1} = \bar D_\ast(u)$, $\bar N_\ast(u) (n^{-1} \log p)^{1/2}
= \bar D_\ast(u)$, and $H_\ast(u) = p u^{1-q/2} n^{(1-q)/2} = \bar
D_\ast(u)$, respectively; let $u_4 = \max(u_1, u_2, u_3, (n^{-1}
\log p)^{1/2})$. For a sufficiently large constant $C_2 > 0$,
$G_\ast(C_2 u_4) = o(D_\ast(u_4))$ and hence the right-hand side
of (\ref{eq:F06626}) is of order $D_\ast(u_4) = O(\tilde M
u_4^{1-r})$ if $u = C_2 u_4$. Let $u_1' = (\tilde M p
n^{1-q})^{1/(q+r)}$ and $u_1'' = (p^{ 1 + 2/q} \tilde M^{-1}
n^{1/q-1})^{1/(1-r)}$. Note that $u_1 = u_1'$ if $p \ge M
(u_1')^{-r}$ and $u_1 = u_1''$ if $p \le M (u_1')^{-r}$. In both
cases we have by elementary calculations that $u_1 = o(u_{\BL})$.
Similarly, we have $u_2 = o(u_{\BL})$ and $u_3 = o(u_{\BL})$.
Hence $u_4 = o(u_{\BL})$ and our rate of convergence $D_\ast(u_4)$
is sharper.

Based on Theorem~\ref{thmm:spectral} and the above discussion, we
have:

\begin{cor}
\label{cor:SpectralF08} Let the conditions in Theorem~\ref{thmm:F08122}
be satisfied and $\alpha> 1/2-1/q$. Let $\Lambda= \sup_{\Sigma
\in{\cal G}_r(\tilde M)} \inf_{u > 0} \|\rho(T_u(\hat{\Sigma}_n)
- \Sigma)\|_2$. Assume $\tilde M \asymp p^\theta$, $0 \le\theta
\le1$ and $p \asymp n^\tau$, $\tau> 0$. Let $\phi_1' = (\tau
\theta+ \tau+ 1 -q) / (q+r)$, $\phi_1'' = (\tau(1-\theta+2/q) -
1 + 1/q) / (1-r)$, $\phi_1 = \min(\phi_1', \phi_1'')$, $\phi_3 =
(2\tau- 2 \tau\theta+ 1-q) / (q-2r)$ and $\phi= \max(\phi_1,
\phi_3)$. \textup{(i)} If $\phi> -1/2$, then $\Lambda=
O(n^{\phi(1-r)+\theta\tau})$. \textup{(ii)} If $\phi\le-1/2$, then
$\Lambda= O(n^{\theta\tau} (n^{-1} \log p)^{(1-r)/2})$.
\end{cor}

\section{Precision matrix estimation for high-dimensional
stationary processes} \label{sec:precision_statproc}

As a straightforward estimate for precision matrices, one can
invert the regularized covariance matrix estimates. However, this
inversion procedure may cause the precision matrix estimate to lose
sparsity. Sparsity of the precision matrix $\Omega= \Sigma^{-1}$
has important statistical meaning because a zero entry in $\Omega
= (\omega_{jk})_{1\le j,k\le p}$ reflects the conditional
independence when $\vz_i$ are multivariate Gaussian. In the
graphical model representation, $\omega_{ij} = 0$ indicates that
there is a missing edge between node~$i$ and node $j$. Performance
bounds for estimating $\Omega$ under dependence is useful for
statistical learning\vadjust{\goodbreak} problems. For direct estimation of precision
matrices that can preserve sparsity, one can adopt entry-wise
1-norm penalized likelihood approaches; see
\citet{friedmanhastietibshirani2008a,banerjeeelghaouidaspremont2008a,ravikumarwainwrightraskuttiyu2008a,rothmanbickellevinazhu2008a,fanfengwu2009a}, which we refer them as Lasso-type precision
matrix estimators. \citet{friedmanhastietibshirani2008a} proposed a
graphical Lasso model and developed a computationally efficient
and scalable algorithm for estimating large precision matrices.
This 1-norm penalized multivariate Gaussian likelihood approach
was also considered by \citet{banerjeeelghaouidaspremont2008a}.
Consistency of the graphical Lasso were studied
in~\citet{rothmanbickellevinazhu2008a,ravikumarwainwrightraskuttiyu2008a}.

The precision matrix estimation procedure considered here is the
graphical Lasso model [\citet{friedmanhastietibshirani2008a}] which
minimizes the objective function
%
\begin{equation}
\label{eq:graphicalLasso} \hat{\Omega}_n(\lambda) = \mathop{\arg\min}
_{\Psi\succ0} \bigl\{\tr(\Psi\hat\Sigma_n) - \log\det(
\Psi) + \lambda|\Psi|_1\bigr\},
\end{equation}
where $\lambda$ is the penalty to be determined later. In
(\ref{eq:graphicalLasso}) $\Psi\succ0$ means that $\Psi$ is
positive-definite. Here we assume the maximum eigenvalue
%
\begin{equation}
\label{eq:uniformclass_Omega} \rho(\Omega) \le\varepsilon_0^{-1} \qquad\mbox{for some } \varepsilon_0 > 0,
\end{equation}
or equivalently the minimum eigenvalue of $\Sigma$ is larger than
$\varepsilon_0$. Note that we do not assume the minimum eigenvalue
of $\Omega$ is uniformly bounded below\vspace*{1pt} from zero. To introduce an
asymptotic theory for the estimate $\hat{\Omega}_n$, we recall
(\ref{eqn:D_u}) and (\ref{eqn:D_u1}) of Theorem~\ref{thmm:F08122}
for the definition of the functions $H(\cdot)$ and $G(\cdot)$ and
also $\tilde\alpha$ and $\tilde\beta$. An analogue of the function
$D(\cdot)$ in this context is
%
\begin{equation}
\label{eq:A08206} D^*(u) = {1\over p^2} \sum
_{j,k=1}^p u \bigl(u \wedge|\omega_{jk}|\bigr).
\end{equation}
Recall Corollary~\ref{cor:M140825} for $\tilde{G}(\cdot)$.

It is interesting and surprising to note that the structure of
Theorem~\ref{thmm:Inv} is very similar to that in Theorem~\ref{thmm:F08122}.
However, the main idea for the proof of Theorem~\ref{thmm:Inv} seems quite different, and our key argument here is
based on convex minimization. It is also interesting to note that
our rate of convergence is expressed in terms of the ${\cal L}^2$
norm; see (\ref{eq:A805948}), while in the previous literature
probabilistic bounds are obtained; see
\citet{ravikumarwainwrightraskuttiyu2008a,rothmanbickellevinazhu2008a,lamfan2009a}. The constant $C$ in
Theorem~\ref{thmm:Inv} can be the same as the one in Theorem~\ref{thmm:F08122}.

\begin{thmm}
\label{thmm:Inv} Let the moment and the dependence conditions in
Theorem~\ref{thmm:F08122} be satisfied and $\lambda= 4 u$. Then
%
\begin{equation}
\label{eq:A805948} \frac{1}{p^2}\E\bigl|\hat\Omega_n(\lambda) -
\Omega\bigr|_F^2 \lesssim D^*(u) + \min \biggl(
{1\over n}, {u^{2-q} \over n^{q/2}}, H(u) + G(C u) \biggr),
\end{equation}
where $C$ is independent of $u, n$ and $p$. Let $u_\flat$ be the
solution to the equation
%
\begin{equation}
\label{eqn:A805928} D^*(u_\flat) = \min\bigl(n^{-1}, \max\bigl(
\tilde{G}(u_\flat), H(u_\flat)\bigr)\bigr).
\end{equation}
Then $\inf_{\lambda> 0} p^{-2}\E|\hat\Omega_n(\lambda) -
\Omega|_F^2 \lesssim D^*(u_\flat)$.
\end{thmm}

\begin{rmk}\label{rem:A8091031}
As an immediate consequence of Theorem~\ref{thmm:Inv}, if the
entries $\omega_{j k}$ of the inverse matrix $\Omega$ satisfy
(\ref{eqn:sparsity_def}) with $0 \le r < 1$, then we have by the
argument in (\ref{eqn:A02849}) that $D^*(u) \le C u^{2-r} M/p^2$.
Similarly, if $\omega_{j k}$ satisfy (\ref{eqn:Expsparsity}), then
$D^*(u) \le C u^2 M \log^r(2+u^{-1})$. Therefore Corollaries
\ref{cor:F_stationary} and \ref{cor:Fexp} are still valid in the
context of precision matrix estimation. 
\end{rmk}

\begin{pf*}{Proof of Theorem \ref{thmm:Inv}}
Using $\Psi= \Omega+ \Delta$, we see that $\hat\Delta_n =
\hat\Omega_n(\lambda) - \Omega$ minimizes
\[
G(\Delta) = \tr(\Delta\hat \Sigma_n)-\log\det(\Psi)+\lambda|
\Psi|_1 + \log\det(\Omega) - \lambda|\Omega|_1.
\]
Hence $G(\hat\Delta_n) \le G(0) = 0$. Let $\Omega_v = \Omega+ v
\Delta$. By Taylor's expansion,
%
\begin{eqnarray}
\label{eqn:G(Delta)} G(\Delta) &= & \tr\bigl[\Delta(\hat \Sigma_n-\Sigma)
\bigr] + \lambda\bigl(|\Omega+\Delta|_1-|\Omega|_1\bigr)
\nonumber
\\[-8pt]
\\[-8pt]
\nonumber
& & {}+
\vec(\Delta)^\top \biggl[\int_0^1(1-v)
\Omega_v^{-1} \otimes\Omega_v^{-1}
\,dv \biggr]\vec(\Delta),
\end{eqnarray}
where $\otimes$ denotes the Kronecker product. Write $\Xi= \hat
\Sigma_n - \Sigma= (\xi_{j k})$, $\mathcal{S}_u = \{(j, k)\dvtx |\omega_{jk}|\ge u\}$ and $\mathcal{W}_u = \{(j, k)\dvtx |\xi_{jk}|\ge
u\}$. Let ${\mathcal{W}^c_u}$ be the complement of
$\mathcal{W}_u$. Then
%
\begin{equation}
\tr(\Delta\Xi) = \tr(\Delta\Xi_{\mathcal{W}_u}) + \tr(\Delta\Xi_{\mathcal{W}^c_u})
\ge- |\Delta|_F | \Xi_{\mathcal{W}_u}|_F - u |
\Delta|_1,
\end{equation}
where the matrix $\Xi_{\mathcal{W}_u} = (\xi_{j k} \mathbf{1}_{(j, k)
\in\mathcal{W}_u})_{1 \le j, k \le p}$. Assume $\alpha> 1/2 -
1/q$. By~(\ref{eq:A8051043}),
%
\begin{equation}
\label{eq:A8071127} \E\bigl(| \Xi_{\mathcal{W}_u}|_F^2\bigr)
\lesssim p^2\bigl( n^{1-q} u^{2-q} +
\bigl(n^{-1}+ u^2\bigr)e^{-C_4n u^2}\bigr) =:
N(u)^2.
\end{equation}
Using the arguments for Theorem~1 in
\citet{rothmanbickellevinazhu2008a}, we have by
(\ref{eq:uniformclass_Omega}) that
%
\begin{equation}
\label{eqn:logdet_term} \vec(\Delta)^\top \biggl[\int_0^1(1-v)
\Omega_v^{-1} \otimes\Omega_v^{-1}
\,dv \biggr]\vec(\Delta) \ge\frac{1}{4}\varepsilon_0^2|
\Delta|_F^2,
\end{equation}
and by letting the penalty $\lambda= 4 u$ that
%
\begin{eqnarray}
\label{eq:A8051132} &&\lambda\bigl(|\Omega+\Delta|_1-|\Omega|_1\bigr) - u
|\Delta|_1 \nonumber\\
&&\qquad\ge\lambda\bigl(\bigl|\Delta^-_{{\cal S}_u^c}\bigr|_1
- 2 |\Omega_{{\cal S}_u^c}|_1 - \bigl|\Delta^+\bigr|_1 - \bigl|
\Delta^-_{{\cal S}_u}\bigr|_1 \bigr) - u |\Delta|_1
\\
&&\qquad
\ge 3 u \bigl|\Delta^-_{{\cal S}_u^c}\bigr|_1 - 8 u |\Omega_{{\cal S}_u^c}|_1
- 5 u \bigl(\bigl|\Delta^+\bigr|_1 + \bigl|\Delta^-_{{\cal S}_u}\bigr|_1
\bigr),\nonumber
\end{eqnarray}
where, for a matrix $\Sigma$, $\Sigma^+ = \diag(\Sigma)$ and
$\Sigma^- = \Sigma- \Sigma^+$. By the Cauchy--Schwarz inequality,
$|\Delta^+|_1 + |\Delta^-_{{\cal S}_u}|_1 \le|\Delta|_F \sqrt{
s_u}$, where $s_u = \# {\cal S}_u$. By
\mbox{(\ref{eqn:G(Delta)})--(\ref{eq:A8051132})},
%
\begin{equation}
G(\Delta) \ge\tfrac{1}{4}\varepsilon_0^2|
\Delta|_F^2 - |\Delta|_F |
\Xi_{\mathcal{W}_u}|_F - 8 u |\Omega_{{\cal S}_u^c}|_1
- 5 u |\Delta|_F \sqrt{s_u}.\vadjust{\goodbreak}
\end{equation}
Since $G(\hat\Delta_n) \le0$, there exists a deterministic
constant $C > 0$ such that
%
\begin{equation}
\label{eq:A8051141} |\hat\Delta_n|_F^2 \le C
\bigl( | \Xi_{\mathcal{W}_u}|_F^2 + u^2
s_u + u |\Omega_{{\cal S}_u^c}|_1\bigr) \le C\bigl(|
\Xi_{\mathcal{W}_u}|_F^2 + p^2 D^*(u)\bigr).
\end{equation}
Then (\ref{eq:A805948}) follows from (\ref{eq:A8071127}) and by
choosing $u$ to minimize the right-hand side of
(\ref{eq:A8051141}); see the argument in (\ref{eqn:Fmax}). The
case with $\alpha\le1/2 - 1/q$ can be similarly handled with
special care (\ref{eq:F24112}) being taken into
(\ref{eq:A8071127}).
\end{pf*}

\citet{ravikumarwainwrightraskuttiyu2008a} studied the graphical
Lasso estimator with off-diagonal entries penalized by the 1-norm.
For i.i.d. $p$-variate vectors with polynomial moment condition,
they showed that if $p = O((n/d^2)^{q/(2\tau)})$ for some
$\tau>2$, where $d$ is the maximum degree in the Gaussian
graphical model, then
%
\begin{equation}
\label{eqn:rate_ravikumar} \frac{1}{p^2}|\hat\Omega_n-\Omega|_F^2
= O_P \biggl( \frac{s+p}{p^2} \cdot\frac{p^{2\tau/q}}{n} \biggr),
\end{equation}
where $s$ is the number of nonzero off-diagonal entries in
$\Omega$. For $\Omega\in\mathcal{H}_0(M)$, we can choose $M =
s+p$. Note that $d \ge s/p$ and thus $d + 1 \ge M/p$. By Remark~\ref{rem:A8091031}, Corollary~\ref{cor:F_stationary} holds. Under
case (ii) [resp., (iii)], our rate of convergence is
$(M/p^2)^{1-2/q} n^{2/q-2}$ [resp., $n^{-1} (\log p) M/p^2$].
Elementary calculations show that both of our rates are of order
$o(M p^{-2} n^{-1} p^{2\tau/q})$. Hence our bounds are much better
than (\ref{eqn:rate_ravikumar}), the one obtained in
\citet{ravikumarwainwrightraskuttiyu2008a}.

We now compare our results with the CLIME (constrained
$L_1$-minimization for inverse matrix estimation) method, a
non-Lasso type estimator proposed in~\citet{cailiuluo2011a}, which
is to
%
\begin{equation}
\label{eqn:CLIME} \mbox{minimize}\quad |\Theta|_1 \quad\mbox{subject to } |\hat
\Sigma_n \Theta-I|_\infty\le \lambda_n.
\end{equation}
\citet{cailiuluo2011a} showed that with $n$ i.i.d. $p$-variate
observations, if $p = o(n^{q/2-1})$, then the rate of convergence
for the CLIME estimator under the normalized Frobenius norm is
$O(\tilde{C}^{4-2r}\tilde{M}(\log{p}/n)^{1-r/2})$, where
$\tilde{C}$ is the upper bound for the matrix $L_1$-norm on the
true precision matrix, and $\tilde{M}$ is in
(\ref{eqn:strongell_qball}). We see that the rates of convergence
under the normalized Frobenius norm are the same for both papers.
This rate of convergence is in general better than those obtained
for the Lasso-type estimators in the polynomial moment case
[\citet{ravikumarwainwrightraskuttiyu2008a}].

\begin{rmk}
Following \citet{rothmanbickellevinazhu2008a}, we can consider the
slightly modified version of the graphical Lasso: let $V = \operatorname{
diag}(\sigma_{1 1}^{1/2}, \ldots, \sigma_{p p}^{1/2})$ and $R$ be
the correlation matrix; let $\hat V$ and $\hat R$ be their sample
versions, respectively. Let $K = R^{-1}$. We estimate $\Omega=
{V}^{-1} K {V}^{-1}$ by $\hat\Omega_\lambda= \hat{V}^{-1}
\hat{K}_\lambda\hat{V}^{-1}$, where
%
\begin{equation}
\label{eqn:graphical_Lasso_variant_step1} \hat{K}_\lambda= \mathop{\arg\min}_{\Psi\succ0} \bigl
\{\tr(\Psi\hat{R}) - \log\det(\Psi) + \lambda\bigl|\Psi^-\bigr|_1 \bigr\}.
\end{equation}
Let $D^-(u) = p^{-2} \sum_{1 \le j \neq k \le p} u(u \wedge
|\omega_{jk}|)$. Using the arguments of Theorem~2 in
\citet{rothmanbickellevinazhu2008a}, we have the following result
on the spectral norm rate of convergence of $\hat\Omega_\lambda$:
Assuming the moment and dependence conditions in Theorem~\ref{thmm:Inv} are satisfied and $\varepsilon_0 \le\rho(\Omega)
\le\varepsilon_0^{-1}$, and then
%
\begin{equation}
\label{eqn:spectral-rate-precision-mat} {
{\rho^2(\hat\Omega_\lambda- \Omega)} \over p^2} \lesssim_{\Prob} D^-(\lambda) +
\min \biggl( {1 \over n}, {\lambda^{2-q} \over
n^{q/2}}, H(\lambda) + G(C
\lambda) \biggr)
\end{equation}
holds if $\max[p^{1/q} n^{-1+1/q}, (\log{p}/n)^{1/2}] \lesssim
\lambda$. Details of the derivation of
(\ref{eqn:spectral-rate-precision-mat}) is given in the
supplementary material [\citet{supp}]. If $\Omega$ satisfies $|\{(j,k) \dvtx \omega_{jk} \neq0, j \neq k\}| \le s$
[\citet{rothmanbickellevinazhu2008a}], we have $\Omega\in\calH_0(
M)$ with $M = s+p$. Simple calculations show that, if $\alpha>
1/2 - 1/q$ and $s = O(p)$, then for $\lambda_\sharp\asymp\max[
(\log{p} / n)^{1/2}, (s^{-1} p^2 n^{1-q})^{1/q} ]$, we have by
(\ref{eqn:spectral-rate-precision-mat}) that
$\rho(\hat\Omega(\lambda_\sharp) - \Omega) = O_{\Prob}(\sqrt{s}
\lambda_\sharp)$, and it reduces to Theorem~2 in
\citet{rothmanbickellevinazhu2008a}.
\end{rmk}

\section{Evolutionary covariance matrix estimation for nonstationary
high-\break dimensional processes}
\label{sec:covariancematrixestiamtion_nonstatproc}

The time series processes considered in Sections~\ref{sec:stationary} and~\ref{sec:precision_statproc} are
stationary. In many situations the stationarity assumption can be
violated, and the graphical structure is time-varying. One may
actually be interested in how the covariance matrices and
dependence structures vary with respect to time.
\citet{zhoulaffertywasserman2010a} and \citet{kolarxing2011a}
studied the estimation of covariance matrices for independent,
locally stationary Gaussian processes. Both requirements can be
quite restrictive in practice.

Here we shall consider nonstationary processes that can be both
dependent and non-Gaussian with mild moment conditions, thus
having a substantially broader spectrum of applicability. To allow
such nonstationary processes, following the framework
in~\citet{draghicescuguillaswu2009a}, we shall consider locally
stationary process
%
\begin{equation}
\label{eq:lsp} \vz_i = \vg(\calF_i; i/n), \qquad 1\le i \le
n,
\end{equation}
where $\vg(\cdot, \cdot) = (g_1(\cdot, \cdot), \ldots, g_p(\cdot,
\cdot))^\top$ is a jointly measurable function such that the
uniform stochastic Lipschitz continuity holds: there exists $C >
0$ for which
%
\begin{equation}
\label{eq:ulc} \max_{j \le p} \bigl\|g_j(
\calF_0; t) - g_j\bigl(\calF_0;
t'\bigr)\bigr\| \le C \bigl|t-t'\bigr| \qquad\mbox{for all } t,
t' \in[0, 1].
\end{equation}
In Examples \ref{ex:A08152}--\ref{ex:A08154} below we present some
popular models of locally stationary processes. Let
$\vz^\diamond_i(t) = \vg(\calF_i; t)$. The preceding condition
(\ref{eq:ulc}) suggests local stationarity in the sense that, for
a fixed $t \in(0, 1)$ and bandwidth $b_n \to0$ with $n b_n \to
\infty$,
%
\begin{equation}
\label{eq:Asplsp} \max_{j \le p} \max_{\lfloor n (t-b_n) \rfloor
\le i \le\lfloor n (t+b_n) \rfloor} \bigl\|
\vz^\diamond_{j, i}(t) - Z_{j, i}\bigr\| \le C b_n
= o(1),
\end{equation}
indicating that the process $(\vz_i)$ over the range $\lfloor n
(t-b_n) \rfloor\le i \le\lfloor n (t+b_n) \rfloor$ can be
approximated by the \textit{stationary} process $\vz^\diamond_i(t)$.
The locally stationarity property suggests that the data
generating mechanism $\vg(\cdot; i/n)$ at time $i$ is close to the
one $\vg(\cdot; i'/n)$ at time $i'$ if $|i-i'|/n$ is small. Hence
the following covariance matrix function is continuous:
%
\begin{equation}
\label{eq:Sigmat} \Sigma(t) = \operatorname{ cov}\bigl(\vg(\calF_0; t)\bigr) =
\E\bigl(\vz(t)\vz(t)^\top\bigr),\qquad t \in(0, 1).
\end{equation}
The covariance matrix $\Sigma_i = \Sigma(i/n)$ of $\vz_i$ can then
be estimated by the approximate stationary process $\vz_l,
\lfloor n (t-b_n) \rfloor\le l \le\lfloor n (t+b_n) \rfloor$, by
using the Nadaraya--Watson or other smoothing techniques. Recall
that in the stationary case the thresholded estimator is defined
as $T_{u}(\hat\Sigma_n) = (\hat\sigma_{jk}
\mathbb{I}(|\hat\sigma_{jk}|\ge u))_{j k}$, where $\hat\Sigma_n =
(\hat\sigma_{jk})$ is the sample covariance matrix given in
(\ref{eq:samplecovariancematrix}). To estimate $\Sigma(t)$, we
substitute $\hat\Sigma_n$ by the kernel smoothed version
%
\begin{equation}
\label{eq:kernelestimator_Sigma}\qquad \hat\Sigma_n(t) = \sum_{m=1}^n
w_m(t) \vz_m\vz_m^\top \qquad\mbox{where } w_m(t) = {{K ({(t-m/n)}/{b_n} )}
\over{\sum_{m=1}^n K ({(t-m/n)}/{b_n} )}}.
\end{equation}
Write $\hat\Sigma_n(t) = (\hat\sigma_{jk}(t))_{j k}$. In
(\ref{eq:kernelestimator_Sigma}), $K(\cdot)$ is a symmetric,
nonnegative kernel with bounded support in $[-1,1]$ and
$\int_{-1}^1 K(v) \,dv=1$. As per convention, we assume that the
bandwidth $b_n$ satisfies the natural condition: $b_n \to0$ and
$n b_n \to\infty$. The thresholded covariance estimator for
nonstationary processes is then defined as
\[
T_{u}\bigl(\hat\Sigma_{n}(t)\bigr) = \bigl(\hat
\sigma_{jk}(t)\mathbb{I}\bigl(\bigl|\hat\sigma_{jk}(t)\bigr|\ge u
\bigr)\bigr)_{1\leq
j,k\leq p}.
\]
Parallelizing Theorem~\ref{thmm:F08122}, we give a general result
for the thresholded estimator for time-varying covariance matrices
of the nonstationary, nonlinear high-dimensional time series. As
in (\ref{eq:functiondependencemeasure_stat}) and
(\ref{eq:srdtail}), we similarly define the functional dependence
measure
%
\begin{equation}
\label{eq:A08141} \theta_{i,w,j} = \max_{0 \le t \le1}
\bigl\|Z_{ji}(t)-Z'_{ji}(t)\bigr\|_w,
\end{equation}
where $Z'_{ji}(t) = g_j(\calF_i', t)$. We also assume that
(\ref{eq:srdtail}) holds. For presentational simplicity let
$\alpha> 1/2-1/q$. Let $n_\sharp= n b_n$, $H_\sharp(u) = {
u^{2-q} n_\sharp^{1-q}}$,
%
\begin{equation}
\label{eqn:A81049} D(u) = {1\over p^2} \max_{0 \le t \le1} \sum
_{j,k=1}^p \bigl(u^2 \wedge
\sigma_{jk}(t)^2\bigr),\qquad G_\sharp(u) =
\bigl(n_\sharp^{-1} + u^2\bigr) e^{-n_\sharp u^2}.
\end{equation}
Theorem~\ref{th:A08946} provides convergence rates for the
thresholded covariance matrix function estimator $T_{u}(\hat
\Sigma_{n}(t))$. Due to the nonstationarity, the bound is worse
than the one in Theorem~\ref{thmm:F08122} since we only use data in
the local window $[n(t-b_n), n(t+b_n)]$. Therefore, in the
nonstationary case a larger sample size is needed for achieving
the same level of estimation accuracy.

\begin{thmm}
\label{th:A08946} Assume $\max_{1 \le j,k \le p} \sup_{t \in
[0,1]} |\sigma''_{jk}(t)| < \infty$ and $\alpha> 1/2 - 1/q$.
Under the moment and dependence conditions of Theorem~\ref{thmm:F08122}, we have
%
\begin{equation}
\label{eq:A808148}
\quad {{\E|T_{u}(\hat\Sigma_{n}(t))-\Sigma(t)|_F^2} \over p^2} \lesssim D(u) + \min\bigl(n_\sharp^{-1},
H_\sharp(u) + G_\sharp(C u)\bigr) + b_n^4
\end{equation}
uniformly over $t \in[b_n, 1-b_n]$, where $C$ is independent of
$u, n, b_n$ and $p$.
\end{thmm}

\begin{pf}
Let $\Sigma^\circ_n(t) = \E\hat\Sigma_n(t) = (\sigma^\circ_{j
k}(t) )_{j k}$. Under the condition on $\sigma''_{jk}(t)$, we have
$\sigma^\circ_{jk}(t) - \sigma_{jk}(t) = O(b_n^2)$ uniformly over
$j, k$ and $t \in[b_n, 1-b_n]$. Hence $|\Sigma^\circ_n(t) -
\Sigma(t)|_F^2 / p^2 = O(b_n^4)$. It remains to deal with $\E
|T_u(\hat\Sigma_n(t)) - \Sigma^\circ_n(t)|_F^2$. With a careful
check of the proof of Theorem~\ref{thmm:F08122}, if we replace
$\hat\sigma_{j k}$ and $\sigma_{j k}$ therein by $\hat
\sigma_{jk}(t)$ and $\sigma^\circ_{jk}(t)$, respectively, then we
can have
%
\begin{equation}
\label{eq:A808119} {{\E|T_{u}(\hat\Sigma_{n}(t))-\Sigma^\circ(t)|_F^2}\over p^2} \lesssim D(u) + \min\bigl(n_\sharp^{-1},
H_\sharp(u) + G_\sharp(C u)\bigr)
\end{equation}
if the following Nagaev inequality holds:
%
\begin{equation}
\label{eq:A08121} \Prob\bigl(\bigl|\hat{\sigma}_{jk}(t)-\sigma^\circ_{jk}(t)\bigr|
> v\bigr) \le\frac{C_2n_\sharp}{(n_\sharp v)^q} + C_3e^{-C_4n_\sharp v^2}.
\end{equation}
The above inequality follows by applying the nonstationary Nagaev
inequality in Section~4 in \citet{liuxiaowu2012a} to the process
$X_m = K ( (t-m/n) / b_n ) (Z_{m j} Z_{m k} - \E(Z_{m j} Z_{m
k}))$, $\lfloor n (t-b_n) \rfloor\le m \le\lfloor n (t+b_n)
\rfloor$. Note that the functional dependence measure of the
latter process is bounded by $\mu(\theta_{i, 2 q, j} + \theta_{i,
2 q, k}) \sup_u |K(u)|$; see (\ref{eq:Julypdm}) and
(\ref{eq:A08141}).
\end{pf}

\begin{rmk}
If in (\ref{eq:kernelestimator_Sigma}) we use the local linear
weights [\citet{MR1383587}], then it is easily seen based on the
proof of Theorem~\ref{th:A08946} that (\ref{eq:A808148}) holds
over the whole interval $t \in[0, 1]$, and the boundary effect is
removed. This applies to the Theorem~\ref{thmm:A08209} below as
well. 
\end{rmk}

A similar result can be obtained for estimating evolutionary
precision matrices of high-dimensional nonstationary processes
$\Omega(t) = \Sigma^{-1}(t)$ where $\Sigma(t)$ is given in
(\ref{eq:Sigmat}). As in the stationary case, we assume that
$\Omega(t)$ satisfies (\ref{eq:uniformclass_Omega}) for all $t \in
[0, 1]$. The actual estimation procedure of $\Omega(t)$ based on
the data $Z_{p \times n}$ is a variant of the graphical Lasso
estimator of $\Omega$, which minimizes the following objective
function:
%
\begin{equation}
\label{eq:A08204} \hat{\Omega}_n(t; \lambda) = \mathop{\arg\min}
_{\Psi\succ0} \bigl\{\tr\bigl(\Psi\hat\Sigma_n(t)\bigr)-
\log\det(\Psi) + \lambda|\Psi|_1\bigr\},
\end{equation}
where $\hat\Sigma_n(t)$ is the kernel smoothed sample covariance
matrix given in (\ref{eq:kernelestimator_Sigma}). The same
minimization program is also used
in~\citet{zhoulaffertywasserman2010a,kolarxing2011a}. As in
(\ref{eq:A08206}) and (\ref{eqn:A81049}), let
%
\begin{equation}
D^*(u) = {1\over p^2} \max_{0 \le t \le1} \sum
_{j,k=1}^p u\bigl(u \wedge\bigl|\omega_{jk}(t)\bigr|
\bigr).
\end{equation}
As in (\ref{eqn:A805928}), choose $\lambda= 4 u_\flat^\sharp$.
For the estimator (\ref{eq:A08204}), we have the following
theorem. We omit the proof since it is similar to the one in
Theorems \ref{thmm:Inv} and \ref{th:A08946}.

\begin{thmm}
\label{thmm:A08209} Assume $\max_{1 \le j,k \le p} \sup_{t \in
[0,1]} |\omega''_{jk}(t)| < \infty$ and $\alpha> 1/2 - 1/q$.
Under the moment and dependence conditions of Theorem~\ref{thmm:F08122}, we have
%
\begin{equation}
\label{eq:A808210} {{\E|\hat\Omega_n(t; 4 u) - \Omega(t)|_F^2}\over p^2} \lesssim D^*(u) + \min
\bigl(n_\sharp^{-1}, H_\sharp(u) + G_\sharp(C
u)\bigr) + b_n^4
\end{equation}
uniformly over $t \in[b_n, 1-b_n]$, where $C$ is independent of
$u, n, b_n$ and $p$. Let $u_\flat^\sharp\ge n_\sharp^{-1/2}$ be
the solution to the equation $\max( G_\sharp(u), H_\sharp(u)) =
D^*(u)$. Then $\inf_{\lambda> 0} p^{-2}\E|\hat\Omega_n(t;
\lambda) - \Omega(t)|_F^2 \lesssim D^*(u_\flat^\sharp)$.
\end{thmm}

\begin{exmp}[{(Modulated nonstationary process [\citet{adak1998a}])}]\label{ex:A08152}
 Let $({\mathbf{y}}_i)$ be a stationary $p$-dimensional process with mean $0$ and
identity covariance matrix. Then the modulated process
%
\begin{equation}
\vz_i = \Sigma^{1/2}(i/n)\mathbf{y}_i,
\end{equation}
has covariance matrix $\Sigma_i = \Sigma(i/n)$.
\citet{zhoulaffertywasserman2010a} considered the special setting
in which $\mathbf{y}_i$ are i.i.d. standard Gaussian vectors, and hence
$\vz_i$ are independent.
\end{exmp}

\begin{exmp}[(Nonstationary linear process)]\label{ex:A08153}
Consider the nonstationary linear
process
%
\begin{equation}
\label{eq:nonstationarylp} \vz_i = \sum_{j=0}^\infty
A_j(i/n) \ve_{i-j}, \qquad 1 \le i \le n,
\end{equation}
where $A_j(\cdot)$ are continuous matrix functions. We can view
(\ref{eq:nonstationarylp}) as a time-varying version of
(\ref{eq:linearprocess}), a framework also adopted
in~\citet{dahlhaus1997a}. As in Example~\ref{exmp:linstat}, we
assume a uniform version
%
\begin{equation}
\max_{k \le p} \sum_{l=1}^p
\max_{0 \le t \le1}a_{j, k l}(t)^2 = O
\bigl(j^{-2-2\gamma}\bigr),\qquad  \gamma> 0.
\end{equation}
\end{exmp}

\begin{exmp}[(Markov chain example revisited: Nonstationary version)]\label{ex:A08154}
 We
consider a nonstationary nonlinear example adapted from Example~\ref{exmp:nonlin_stat}. Let the process $(\vz_i)$ be defined by
the iterated random function
%
\begin{equation}
\label{eq:iteratedrandomfun_nonstat} \vz_i=\vg_i(\vz_{i-1},
\ve_i),
\end{equation}
where $\vg_i(\cdot,\cdot)$ is an $\mathbb{R}^p$-valued and jointly
measurable function that may change over time. As in Example~\ref{exmp:nonlin_stat}, we assume $\vg_i$ satisfy: (i) there
exists some $\vx_0$ such that $\sup_i \| \vg_i(\vx_0,\ve_0)\|_{2q}
< \infty$; (ii)
\[
L:=\sup_i \E|L_i|^q < 1 \qquad\mbox{where } L_i = \sup_{\vx\neq
\vx'}\frac{\|\vg_i(\vx,\ve_0)-\vg_i(\vx',\ve_0)\|_{2q}}{|\vx
-\vx'|}.\vadjust{\goodbreak}
\]
Then $(\vz_i)$ have the GMC property with $\Theta_{m,2q} =
O(L^m)$. Therefore, Theorem~\ref{th:A08946} can be applied with
$\alpha> 1/2-1/q$ and $\tilde\beta= 1$.
\end{exmp}

\section*{Acknowledgments}

We thank two anonymous referees, an  Associate Editor and the Editor for their helpful
comments that have improved the paper.

\begin{supplement}[id=suppA]
\stitle{Additional proofs}
\slink[doi]{10.1214/13-AOS1182SUPP} 
\sdatatype{.pdf}
\sfilename{aos1182\_supp.pdf}
\sdescription{The supplementary file contains the proof of
relation (\ref{eqn:spectral-rate-precision-mat}):
spectral norm convergence rate for precision matrix.}
\end{supplement}

%


\printaddresses

\end{document}